\documentclass[12pt]{amsart}

\usepackage{amscd}
\usepackage{pb-diagram}

\newtheorem{thm}{Theorem}[section]   % Numbered within each section
\newtheorem{cor}[thm]{Corollary}     % Numbered along with thm
\newtheorem{lem}[thm]{Lemma}         % Numbered along with thm
\newtheorem{prop}[thm]{Proposition}  % Numbered along with thm

\theoremstyle{definition}

\theoremstyle{remark}

\newtheorem{rem}[thm]{Remark}        % Numbered along with thm
\numberwithin{equation}{section}     % Number eqns within sections

\newcommand{\secref}[1]{Section~\ref{#1}}
\newcommand{\thmref}[1]{Theorem~\ref{#1}}
\newcommand{\corref}[1]{Corollary~\ref{#1}}
\newcommand{\lemref}[1]{Lemma~\ref{#1}}
\newcommand{\propref}[1]{Proposition~\ref{#1}}

\newcommand{\ind}{\operatorname{Ind}}
\newcommand{\Ind}{\ind}
\newcommand{\res}{\operatorname{Res}}
\newcommand{\Res}{\res}
\newcommand{\Ex}{\operatorname{Ex}}
\newcommand{\Sub}{\operatorname{Sub}}
\newcommand{\rep}{\operatorname{Rep}}
\newcommand{\Rep}{\rep}

\newcommand{\id}{\operatorname{\rm id}}

\newcommand{\im}{\operatorname{\rm im}}

\newcommand{\Ad}{\Ad}

\newcommand{\dashind}{\operatorname{\!-Ind}}
\renewcommand{\inf}{\operatorname{Inf}}

\newcommand{\contains}{\supset}
\renewcommand{\phi}{\varphi}

\newcommand{\half}{\frac{1}{2}}

\renewcommand{\H}{\mathcal H}

\newcommand{\I}{\mathcal I}
\newcommand{\K}{\mathcal K}
\newcommand{\D}{\mathcal D}
\renewcommand{\L}{\mathcal L}

\newcommand{\cstar}{\ensuremath{C^*}-}

\newcommand{\righttext}[1]{\qquad\text{#1 }}

\newcommand{\deltahat}{{\hat{\delta}}}

\newcommand{\lip}[3]{
  {\vphantom\langle}_{#1}\!\!\left\langle{#2},{#3}\right\rangle  }
\newcommand{\rip}[3]{
  \left\langle{#2},{#3}\right\rangle_{\!{#1}}  }

\begin{document}

\title{Duality of restriction and induction for \cstar coactions}

\author{S. Kaliszewski}
\address{Department of Mathematics\\University of Newcastle\\
 Newcastle, New South Wales 2308\\Australia}
\email{kaz@frey.newcastle.edu.au}

\author{John Quigg}
\address{Department of Mathematics\\Arizona State University\\
 Tempe, Arizona 85287}
\curraddr{Department of Mathematics\\University of Newcastle\\
 Newcastle, New South Wales 2308\\Australia}
\email{quigg@math.la.asu.edu}

\author{Iain Raeburn}
\address{Department of Mathematics\\University of Newcastle\\
 Newcastle, New South Wales 2308\\Australia}
\email{iain@frey.newcastle.edu.au}

\thanks{This research was partially supported by the National Science
Foundation under Grant No. DMS9401253, and by the Australian Research
Council.}



\date{March 22, 1996 (revised)}

\begin{abstract}
Consider a coaction $\delta$ of a locally compact group $G$ on a
\cstar algebra $A$, and a closed
normal subgroup $N$ of $G$.  We prove, following
results of Echterhoff for abelian $G$, that Mansfield's imprimitivity
between $A\times_{\delta|}G/N$ and $A\times_\delta G\times_{\deltahat,r}N$
implements equivalences between Mansfield induction of representations
from $A\times G/N$ to $A\times G$ and restriction of representations
from $A\times G\times_r N$ to $A\times G$, and between restriction of
representations from $A\times G$ to $A\times G/N$ and Green induction
of representations from $A\times G$ to $A\times G\times_r N$.  This
allows us to deduce properties of Mansfield induction from the known
theory of ordinary crossed products.
\end{abstract}

\maketitle

%====================================================================

\section{Introduction}\label{intro-sec}

In applications of duality theory for crossed products one often
has to know how common constructions such as induction and
restriction of representations behave under duality.  Thus, after
piecemeal results by several authors, Echterhoff was led to prove that,
for systems involving abelian groups, induction and restriction are
dual to one another \cite{EchDI}.  He later used these results to great
effect in his analysis of crossed products with continuous trace
\cite{EchCP}.  

To state Echterhoff's theorem precisely, we fix an
action $\alpha$ of a locally compact abelian group $G$ on a \cstar algebra
$A$, and a closed subgroup $N$ of $G$.  Green's imprimitivity
theorem gives an imprimitivity bimodule $X_N^G$ implementing a Morita
equivalence between $A\times_\alpha N$ and the imprimitivity algebra
$(A\otimes C_0(G/N))\times_{\alpha\otimes\tau}G$; the latter algebra is
naturally isomorphic to the crossed product $A\times_\alpha
G\times_{\hat{\alpha}} N^\perp$ by the dual action $\hat{\alpha}$ of
the subgroup $(G/N)\widehat{\ } = N^\perp$ of $\hat{G}$.  Let $\Ind\colon
\Rep A\times_\alpha N\to \Rep A\times_\alpha G$ denote the map on
equivalence classes of representations given by induction of
representations in the sense of Green: if $\pi$ is a representation of
$A\times_\alpha N$ on $\H$, then $\Ind\pi$ is by definition the natural
left action on (a completion of) $X_N^G\otimes_{A\times_\alpha N}\H$.
Similarly, let $\Res_G^N\colon \Rep A\times_\alpha G\to \Rep
A\times_\alpha N$ denote the map given by restricting a covariant
representation $(\pi,U)$ of $(A,G,\alpha)$ to  the representation
$(\pi,U|_N)$ of $(A,N,\alpha)$.  Then Echterhoff's theorem says that we
have commutative diagrams
\begin{equation}\label{ech-res-ind-diag}
\begin{CD}
\Rep A\times_\alpha G @= \Rep A\times_\alpha G\\
@V{\Res_G^N}VV @VV{\Ind_{\{e\}}^{H^\perp}}V\\
\Rep A\times_\alpha N @>>{X_N^G}> \Rep A\times_\alpha
G\times_{\hat{\alpha}} N^\perp
\end{CD}
\end{equation}
and
\begin{equation}\label{ech-ind-res-diag}
\begin{CD}
\Rep A\times_\alpha G @= \Rep A\times_\alpha G\\
@A{\Ind_N^G}AA @AA{\Res^{\{e\}}_{H^\perp}}A\\
\Rep A\times_\alpha N @>>{X_N^G}> \Rep A\times_\alpha
G\times_{\hat{\alpha}} N^\perp,
\end{CD}
\end{equation}
in which the horizontal arrows on the bottom are the bijections 
induced by the imprimitivity bimodule $X_N^G$.  

As it stands, this theorem only makes sense for abelian groups, because
otherwise there is no dual action.  However, as Echterhoff himself
observed, this is an obvious case of a theorem about abelian groups
which should extend to an arbitrary locally compact group $G$ and a
closed
normal subgroup $N$, replacing dual actions  
by dual coactions.  
(Another such theorem is \cite[Theorem~2.4]{OP-PI}, which was extended to
non-abelian groups in \cite[Theorem~4.4]{QR-IC}.)
This project has been carried
through by Echterhoff and the first and third authors
\cite{EKR}, relating induction and
restriction for $A\times_\alpha G$ to,
respectively, the restriction and induction
processes 
of Mansfield \cite{ManJF}\ 
for crossed products by coactions.

Here we study the dual situation in which we start with a
coaction $\delta$ of a locally compact group $G$ on $A$.
Our main theorems give (in special cases) commutative diagrams
\begin{equation}\label{our-res-ind-diag}
\begin{CD}
\Rep A\times_\delta G @= \Rep A\times_\delta G\\
@V{\Res}VV @VV{\Ind}V\\
\Rep A\times_{\delta|}G/N @>>{Y_{G/N}^G}> \Rep A\times_\delta
G\times_{\hat{\delta},r}N
\end{CD}
\end{equation}
and
\begin{equation}\label{our-ind-res-diag}
\begin{CD}
\Rep A\times_\delta G @= \Rep A\times_\delta G\\
@A{\Ind}AA @AA{\Res}A\\
\Rep A\times_{\delta|}G/N @>>{Y_{G/N}^G}> \Rep A\times_\delta
G\times_{\hat{\delta},r}N
\end{CD}
\end{equation}
similar to \eqref{ech-res-ind-diag}\ 
and \eqref{ech-ind-res-diag}.  This time, the horizontal arrows
are the bijections given by Mansfield's Morita equivalence of
$A\times_{\delta|} G/N$ with $A\times_\delta G\times_{\deltahat,r} N$.
These theorems are potentially very useful, since they reduce questions
about Mansfield induction and restriction to questions about the much
more familiar and manageable theory of Green \cite{GreLS}.
Indeed, in the special case $N=G$, diagrams 
\eqref{our-res-ind-diag}\ and \eqref{our-ind-res-diag}\ 
reduce to results used by Gootman and Lazar \cite{GL-AN}\ to study
induced ideals and the primitive ideal space of crossed products by
coactions. 

We begin with a preliminary section in which we discuss our
conventions about Hilbert modules, basic facts about induction of
representations, and a version of Green's imprimitivity theorem for
reduced crossed products.  This is necessary because we want to avoid
unnecessary assumptions of amenability; for the same reason, we shall
use the version of Mansfield induction given in \cite{KQ-IC}.  We
stress, however, that for amenable groups, the
only part of \cite{KQ-IC}\ needed here is the
extension of Mansfield's theorem to full coactions and twisted systems.  
A special case of this extension
is in \cite{PR-TC}, and the general case is hinted at in
\cite{ER-ST}. 

Our first main results in 
\secref{res-ind-sec}\  concern the ``$\Res$-$\Ind$'' diagram
\eqref{our-res-ind-diag}.  The result is more general than that stated
above in three respects.  First, 
we consider two normal subgroups $N\subset H\subset G$.  Second, we
consider also the twisted crossed products of \cite{PR-CO}, which
introduces another normal subgroup $K$ containing $H$; our strategy,
both here and in \secref{ind-res-sec}, 
is to prove the untwisted version, and then
check that the constructions respect the twists.
Third, we prove that the diagram commutes in
a very strong sense:  all the maps in the diagram are implemented by
Hilbert modules, so their compositions are given by tensor products
of modules, and we prove that the tensor products corresponding to
the two alternative routes are naturally isomorphic.  This not only
gives a stronger theorem
than the corresponding result of \cite{EchDI}, 
but also a distinctly shorter proof.
Thus we can claim that this is another
example of a theorem about abelian groups which becomes much cleaner
when seen from the point of view of coactions.  (Well, we would,
wouldn't we.) 

In \secref{ind-res-sec}\ we turn to the ``$\Ind$-$\Res$'' diagram
\eqref{our-ind-res-diag}.  As in \secref{res-ind-sec}, 
our main results are much
more general than stated above.  
The overall pattern of this section is like
that of \secref{res-ind-sec}, 
but this time passing to twisted crossed products is harder.  
In Echterhoff's paper, the second diagram \eqref{ech-ind-res-diag}\
amounts to the induction in stages of Green \cite{GreLS}.  Mansfield
did not prove such a theorem for his induction process, so we are
forced to do it here 
(\corref{mansfield-stages-thm}).  

In \secref{me-infl-stab-sec}\ 
we give some applications of our theorems.  We
prove that restriction and induction are compatible with Morita
equivalence in general, and in particular with the
stabilization trick of \cite{ER-ST}, which allows us to replace twisted
crossed products by ordinary ones.  (This result for abelian groups was
a fundamental tool in \cite{EchCP}.)  Along the way we deduce, from the
corresponding properties of $\Ind$ and $\Res$ for crossed products by
actions, that both $\Res$ and $\Ind$ are respected by Morita
equivalence of coactions.  
We feel that this is a good illustration of how our results might be
used to deduce information about coactions from known properties of
actions. 

Finally, in the last section, we use our main theorems to
study the maps $\Res$, $\Ind$, $\Ex$ and $\Sub$ on ideals of crossed
products by coactions, obtaining generalizations to the case of
nonamenable $G$ of many of the results in \cite[\S3]{GL-AN}.  Here we
see the power not only of $\Res$-$\Ind$ (and $\Sub$-$\Ex$) duality, but
also that of the Hilbert module techniques, since our maps on spaces of
representations automatically give rise, via the Rieffel
correspondence, to maps on spaces of ideals.  

This research was carried out while the second author was visiting the
University of Newcastle in 1994 and
1995, and while the first author was visiting 
Arizona State University in 1995.  
The various visitors are grateful to their
respective hosts for their hospitality.

%====================================================================

\section{Preliminaries}\label{setup}\label{prelim-sec}

Throughout, $G$ will be a locally compact group with modular function
$\Delta_G$.  We use left Haar
measure.  The group \cstar algebra of $G$
is denoted $C^*(G)$; a subscript $r$, as in $C^*_r(G)$ or $B\times_r
G$, always indicates a reduced object.  
Nondegenerate
homomorphisms of \cstar algebras extend to homomorphisms of their
multiplier algebras, and this will be done implicitly.

\subsection*{Coactions and Imprimitivity}

We use the conventions of 
\cite{KQ-IC}, \cite{QuiFR}, 
\cite[Section~7]{QR-IC}, and \cite{RaeCR}, although the
latter uses maximal tensor products. Our coactions use minimal tensor
products, are injective, and are full, i.e., use $C^*(G)$. Let
$(A,G,\delta)$ be a coaction. We let $(A\times G,j_A,j_G)$ denote the
crossed product, $\hat\delta$ the dual action of $G$ on $A\times G$,
and $\delta^{\rm n}$ the normalization of $\delta$. 
If $N$ is a closed normal subgroup of $G$, we let $\delta|$ denote the
restricted coaction of $G/N$ on $A$ and $A\times G/N$ the restricted
crossed product. Moreover we let $\mu|$ denote the restriction to
$C_0(G/N)$ of a nondegenerate homomorphism $\mu$ of $C_0(G)$.

In \cite{KQ-IC}\ 
the first two authors generalized Mansfield's imprimitivity machine
\cite{ManJF}.
When $(A,G,\delta)$ is a nondegenerate
coaction and $N$ is a closed normal subgroup of $G$, there are dense
$*$-subalgebras $\D$ and $\D_N$ of $A\times G$ and $\im j_A\times j_G|\subset
M(A\times G)$, respectively, such that $\D$ is a (full) pre-Hilbert
$\D_N$-module under right multiplication and inner product
$$\langle x,y\rangle_{\D_N}=\int_N \hat\delta_n(x^*y)\,dn
\righttext{for}x,y\in\D,$$
the integral converging strictly in $M(A\times G)$. Let $Y_{G/N}^G$
denote the completion of the pre-Hilbert module $\D$. Then $Y_{G/N}^G$
has a left $A\times G\times_r N$-module action determined by left
multiplication on $\D$ and
$$n\cdot x=\Delta(n)^{\half}\hat\delta_n(x)
\righttext{for}n\in N,x\in\D.$$
When $N$ is amenable and $\delta$ is a reduced coaction, Mansfield
\cite[Theorem~27]{ManJF}\ proves that $A\times G\times N$ is Morita
equivalent to $A\times G/N$.  For nonamenable subgroups and full
coactions, the corresponding result is:

\begin{thm}\cite[Corollary~3.4]{KQ-IC}\label{manimp}
Let $(A,G,\delta)$ be a nondegenerate
coaction and $N$ a closed normal subgroup of $G$ such that
$$j_A\times j_G|\colon A\times G/N\to M(A\times G)$$
is faithful.
Then $Y_{G/N}^G$ is
an $A\times G\times_r N$ -- $A\times G/N$ imprimitivity bimodule.
\end{thm}

In view of the above theorem, if $\delta$ is a nondegenerate
coaction and $N$ is a closed normal subgroup of $G$, we say
\emph{Mansfield imprimitivity works} for $N$ and $\delta$ 
whenever 
$j_A\times j_G|\colon A\times G/N\to M(A\times G)$ is faithful
\cite[Definition~3.5]{KQ-IC}.
When Mansfield imprimitivity works we let $\langle
\cdot,\cdot\rangle_{A\times G/N}$ denote the extension to $Y_{G/N}^G$ of
the inner product $\langle \cdot,\cdot\rangle_{\D_N}$ on $\D$.
Mansfield's computations show that the left inner product ${}_{A\times
G\times_r N}\langle x,y\rangle$ for $x,y\in\D$ can be identified with
the element
$${}_{A\times G\times_r N}\langle x,y\rangle(n)
=x\hat\delta_n(y^*)\Delta(n)^{-\half}$$
of $C_c(N,\D)$.
When
$\delta$ is nondegenerate, Mansfield imprimitivity works if either
$N$ is amenable or $\delta$ is normal \cite[Lemma~3.2]{KQ-IC}. 
If Mansfield imprimitivity works for
$N$ and $\delta$, then it also works for any closed subgroup contained in $N$
\cite[Theorem~5.2]{KQ-IC}\ and 
any coaction Morita equivalent to $\delta$ \cite[Theorem~5.3]{KQ-IC}.
If $K$ is a closed normal subgroup of $G$ 
containing $N$ and $(A,K,\epsilon)$ is a coaction,
then Mansfield imprimitivity works for $N$ and $\epsilon$ if and only
if it works for $N$ and 
the inflated coaction
$(A,G,\inf\epsilon)$ \cite[Theorem~5.4]{KQ-IC}.

When we say
$(A,G,G/K,\delta,\tau)$ is a twisted coaction, we 
mean $K$ is a closed normal
subgroup of $G$ and $\tau\colon C_0(G/K)\to M(A)$ is a twist for
$\delta$ over $G/K$ \cite{PR-TC}. 
We let $I_\tau$ denote the twisting ideal of
$A\times G$,
$A\times_{G/K}G=(A\times G)/I_\tau$ the twisted crossed product, and
$\tilde\delta$ the dual action of $K$ on $A\times_{G/K}G$. 
If further $N$ is a closed normal
subgroup of $G$ contained in $K$, then $\tau$ is also a twist for the
restricted coaction $(A,G/N,\delta|)$ over the quotient $G/K\cong
(G/N)/(K/N)$. Let $I_\tau^N$ denote the twisting ideal of $A\times
G/N$. Then there is a restricted twisted crossed product
$A\times_{G/K}G/N=(A\times G/N)/I_\tau^N.$
The following result generalizes \cite[Theorem~4.1]{PR-TC}:

\begin{thm}\cite[Theorem~4.4]{KQ-IC}\label{mprimp}
Let $(A,G,G/K,\delta,\tau)$ be a nondegenerate twisted
coaction and $N$ a closed normal subgroup of $G$ contained in $K$ such that
Mansfield imprimitivity works for $N$ and $\delta$. 
Then the quotient
$Z_{G/N}^G=Y_{G/N}^G/(Y_{G/N}^G\cdot I_\tau^N)$ is an
$A\times_{G/K}G\times_r N$ -- $A\times_{G/K}G/N$ imprimitivity bimodule.
\end{thm}

If $(A,G,G/K,\delta,\tau)$ is a nondegenerate twisted coaction, then
Mansfield imprimitivity works for $K$ and $\delta$ if and only if
$\delta$ is normal \cite[Lemma~3.6]{KQ-IC}.
If $(A,G,G/K,\delta,\tau)$ is a nondegenerate normal twisted coaction
and $N$ is a closed normal subgroup of $G$ contained in $K$, then Mansfield
imprimitivity works for $N$ and $\delta$
since $\delta$ is normal, and also for $N$
and the Morita equivalent stabilized coaction 
$(A\times_{G/K}G\times_r K,K,(\widehat{\widetilde\delta})^{\rm n})$ 
\cite[Theorem~3.1]{ER-ST}, \cite[Theorem~5.5]{KQ-IC}.

\subsection*{Hilbert Modules and Rieffel Induction}

Everything in this paper revolves around Rieffel's induction process, so
we should make our conventions explicit. 
For more detailed treatments of this material we refer the reader to
\cite{LanHC}, \cite{RieIR}, \cite{KQ-IC}.
All our Hilbert modules
will be full, i.e., the closed span of the inner product generates the
\cstar algebra. If $X$ is a right Hilbert $B$-module and $A$ acts
nondegenerately on $X$ by adjointable $B$-module maps (so there is a
homomorphism $A\to\L_B(X)$ such that $AX=X$),
we say $X$ is a right-Hilbert $A$ -- $B$ bimodule.
(This terminology first appears in \cite{BuiFC}.)
If $X$ is also a left Hilbert $A$-module
such that $\lip{A}{x}{y}\cdot z = x\cdot \rip{B}{y}{z}$ for $x,y,z\in
X$, then of course $X$ is an
$A$ -- $B$ imprimitivity bimodule. 
We denote the reverse bimodule by $\tilde{X}$, with elements
$\tilde{x}$.  

When $X$ is a right-Hilbert
$A$ -- $B$ bimodule, Rieffel induction gives a functor
$$X\dashind^A_B\colon\rep B\to\rep A,$$
and we leave out parts of the notation if confusion seems unlikely.
Actually, $X\dashind^A_B$ can be factored as
$$X\dashind\colon\rep B\to\rep\K_B(X)$$
followed by ``restriction'' from $\rep\K_B(X)$ to $\rep A$. Since we will
need it a lot, we abstract this latter bit: if $\pi\colon A\to M(B)$ is
a nondegenerate homomorphism, composition with $\pi$ gives a
``restriction'' map
$$\res_B^A = \pi^*\colon\rep B\to\rep A.$$
We view this as a Rieffel induction process: $B$ becomes a 
right-Hilbert $A$ -- $B$ bimodule via
\begin{align*}
\langle b,c\rangle_B&=b^*c\\
b\cdot c&=bc\\
a\cdot b&=\pi(a)b,
\end{align*}
for $a\in A$, $b,c\in B$. 

Now suppose $\pi\colon A\to M(B)$ is a
nondegenerate homomorphism and $X$ is a right-Hilbert
$B$ -- $C$ bimodule. 
Then $X$ becomes a  right-Hilbert $A$ -- $C$ bimodule via
$$a\cdot x=\pi(a)x\righttext{for}a\in A,x\in X.$$
On the other hand, we
can regard $B$ as a  right-Hilbert
$A$ -- $B$ bimodule, 
and the map $b\otimes x\mapsto bx$ induces an isomorphism
$$B\otimes_B X\to X$$
of  right-Hilbert $A$ -- $C$ bimodules. In particular, we get a
commutative diagram 
\begin{equation}\label{restrict}
%----------------------------------------------------------------------
\begin{diagram}
	\node{\rep C}
		\arrow{e,t}{{}_BX_C\dashind}
		\arrow{se,b}{{}_AX_C\dashind}
	\node{\rep B}
		\arrow{s,r}{\Res}\\
	\node[2]{\rep A}
\end{diagram}
%----------------------------------------------------------------------
\end{equation}
where we use prescripts and postscripts to indicate the coefficient
algebras when necessary. In general, we will omit parts of the
notation, so that when we say
\begin{equation}\label{commute}
%----------------------------------------------------------------------
\begin{diagram}
	\node{\rep C}
		\arrow{e,t}{Y}
		\arrow{se,b}{Z}
	\node{\rep B}
		\arrow{s,r}{X}\\
	\node[2]{\rep A}
\end{diagram}
%----------------------------------------------------------------------
\end{equation}
is a commutative diagram, we mean $X$ is a 
right-Hilbert $A$ -- $B$ bimodule, and similarly for $Y$ and $Z$, and 
the equation $(X\dashind)\circ (Y\dashind) = Z\dashind$ holds 
in the strong sense that
$$X\otimes_B Y\cong Z$$
as  right-Hilbert $A$ -- $C$ bimodules.

Recall that Rieffel induction gives rise to maps between ideals, so if
$I$ is an ideal of $B$ and if $\pi$ is any nondegenerate representation
of $B$ with kernel $I$, then $X\dashind I$ 
is the kernel of $X\dashind\pi$. When
we have a commutative diagram of Hilbert modules as in \eqref{commute},
we of course get
$$(X\dashind)\circ(Y\dashind)=Z\dashind$$
as maps from ideals of $C$ to ideals of $A$.

We will often want to pass from a commutative diagram of Hilbert modules to
quotients. There is a subtle point that needs checking:

\begin{lem}
Let $X$ and $Y$ be  right-Hilbert $A$ -- $B$ and  
right-Hilbert $B$ -- $C$ bimodules, 
respectively, and let $K$ be an ideal of $C$.
Further, let $J=Y\dashind K$ and $I=X\dashind J$ be the corresponding
induced ideals of $B$ and $A$, respectively, and let $q_X\colon X\to
X/(X\cdot J)$ and $q_Y\colon Y\to Y/(Y\cdot K)$ be the quotient maps. Then
$q_X\otimes q_Y$ induces an isomorphism between the  
right-Hilbert $A/I$ -- $C/K$ bimodules 
$(X\otimes_B Y)/\left((X\otimes_B Y)\cdot K\right)$ and
$(X/(X\cdot J))\otimes_{B/J} (Y/(Y\cdot K))$.
\end{lem}

\begin{proof}
Straightforward; a slightly different version was given in
\cite[Lemma 1.10]{RaePG}.
\end{proof}

\begin{cor}\label{quotient}
Suppose the diagram
%----------------------------------------------------------------------
\begin{equation}
\begin{diagram}
	\node{\rep C}
		\arrow{e,t}{Y}
		\arrow{se,b}{Z}
	\node{\rep B}
		\arrow{s,r}{X}\\
	\node[2]{\rep A}
\end{diagram}
\end{equation}
%----------------------------------------------------------------------
commutes in the usual strong sense that $Z\cong X\otimes_B Y$.
Further suppose that $K$ is an ideal of $C$, and set $J=Y\dashind K$ and
$I=X\dashind J$. Then the diagram
%----------------------------------------------------------------------
\begin{equation}
\begin{diagram}
	\node{\rep C/K}
		\arrow{e,t}{Y/(Y\cdot K)}
		\arrow{se,b}{Z/Z\cdot K}
	\node{\rep B/J}
		\arrow{s,r}{X/X\cdot J}\\
	\node[2]{\rep A/I}
\end{diagram}
\end{equation}
%----------------------------------------------------------------------
also commutes in the usual strong sense.
\end{cor}

\begin{proof}
Since $Z\cong X\otimes_B Y$, certainly
$$Z/(Z\cdot K)\cong (X\otimes_B Y)/((X\otimes_B Y)\cdot K),$$
so the above lemma immediately gives the corollary.
\end{proof}
Similarly, commutative diagrams with any number of vertices pass to
quotients.  

\subsection*{Green Induction for Reduced Crossed Products}

We will need an induction process for reduced crossed products by
actions.  We could 
deduce it from Green's version by applying 
\cite{QS-RH}, but we give a direct argument
since we need the explicit imprimitivity bimodule. 

Let $(B,G,\alpha)$ be an action and $H$ a closed subgroup of $G$.
Recall that Green's inducing process starts with the 
right-pre-Hilbert 
$C_c(G,B)$ -- $C_c(H,B)$ bimodule $C_c(G,B)$, where the
operations for $f,x,y\in C_c(G,B)$, $g\in C_c(H,B)$ are given by
\begin{align*}
(f\cdot x)(t)
&=\int_G f(s)\alpha_s(x(s^{-1}t))\Delta_G(s)^{\half}\,ds\\
(x\cdot g)(t)
&=\int_H x(th^{-1})\alpha_{th^{-1}}(g(h))\Delta_H(h)^{-\half}\,dh\\
\langle x,y\rangle_{C_c(H,B)}(h)
&=\int_G
\alpha_s\bigl(x(s^{-1})^*y(s^{-1}h)\bigr)\Delta_H(h)^{-\half}\,ds.
\end{align*}
The particular arrangement of modular functions comes from
\cite{RaeIC}; Green uses unorthodox conventions.
Let $Z^G_H$ be the completion of the pre-Hilbert module
$C_c(G,B)$, so $Z^G_H$ is a right-Hilbert
$B\times G$ -- $B\times H$ bimodule. 
We need to know that the kernels of the regular
representations match up, so that $Z^G_H$ passes to a 
right-Hilbert $B\times_r G$ -- $B\times_r H$ bimodule.

\begin{lem}\label{reduced}
Let $(B,G,\alpha)$ be an action, $H$ a 
closed subgroup of $G$, $Z^G_H$ Green's
right-Hilbert 
$B\times G$ -- $B\times H$ bimodule, and $I$ the kernel
of the regular representation of $B\times H$. Then the induced ideal
$Z^G_H\dashind_{B\times H}^{B\times G}I$ is the kernel of the regular
representation of
$B\times G$.
Consequently,  the quotient $X^G_H=Z^G_H/(Z^G_H\cdot I)$ is a 
right-Hilbert $B\times_r G$ -- $B\times_r H$ bimodule.
\end{lem}

\begin{proof}
This follows from induction in stages: 
the proof of \cite[Proposition~8]{GreLS} shows that the diagram
\begin{equation}\label{stages}
%----------------------------------------------------------------------
\begin{diagram}
	\node{\rep B}
		\arrow{e,t}{Z_{\{e\}}^H}
		\arrow{se,b}{Z_{\{e\}}^G}
	\node{\rep B\times H}
		\arrow{s,r}{Z_H^G}\\
	\node[2]{\rep B\times G}
\end{diagram}
%----------------------------------------------------------------------
\end{equation}
commutes in the usual strong sense.
If $\pi$ is any
faithful representation of $B$, then $Z^H_{\{e\}}\dashind\pi$ and
$Z^G_{\{e\}}\dashind\pi$ are the regular representations of $B\times H$
and
$B\times G$, respectively. So the ideals $Z^H_{\{e\}}\dashind\{0\}$
and
$Z^G_{\{e\}}\dashind\{0\}$ are the kernels of the regular
representations.
By commutativity of diagram \eqref{stages},
$$Z^G_{\{e\}}\dashind\{0\}=Z^G_H\dashind
\left( Z^H_{\{e\}}\dashind\{0\}\right),$$
so indeed $Z^G_H$ induces the kernel $I$ of $B\times H\to B\times_r H$
to the kernel of $B\times G\to B\times_r G$.
\end{proof}

\begin{cor}\label{reduced-stages}
Let $(B,G,\alpha)$ be an action, $H$ a closed
subgroup of $G$, and $X^G_H$ the
right-Hilbert $B\times_r G$ -- $B\times_r H$ bimodule obtained in
\lemref{reduced}. Then the diagram
\begin{equation}
%----------------------------------------------------------------------
\begin{diagram}
	\node{\rep B}
		\arrow{e,t}{X_{\{e\}}^H}
		\arrow{se,b}{X_{\{e\}}^G}
	\node{\rep B\times_r H}
		\arrow{s,r}{X_H^G}\\
	\node[2]{\rep B\times_r G}
\end{diagram}
%----------------------------------------------------------------------
\end{equation}
commutes in the usual strong sense.
\end{cor}

\begin{proof}
This follows from the above lemma and \corref{quotient}.
\end{proof}
We emphasize that $X_H^G$ may be viewed as the completion of $C_c(G,B)$
with respect to the norm induced by the $B\times_r H$-valued pre-inner
product.  In particular, the actions of $B\times_r H$ and $B\times_r G$
on $X_H^G$ are determined by the covariant representations of
$(B,H,\alpha)$ and $(B,G,\alpha)$ on $C_c(G,B)$. 

%=====================================================================

\section{Mansfield restriction and Green induction}\label{res-ind-sec}

Suppose we have a
twisted coaction $(A,G,G/K,\delta,W)$, and 
closed normal subgroups $N\subset H$ of $G$ contained in $K$.
In this section, we show that when Mansfield imprimitivity works, 
the following diagram commutes in the
usual strong sense:
%--------------------------------------------------------------------
\begin{equation*}
\begin{CD}
\Rep A\times_{G/K}G/N  	@>>>	\Rep A\times_{G/K}G\times_r N  \\
@V\Res VV			@VV\Ind V\\
\Rep A\times_{G/K}G/H   @>>>    \Rep A\times_{G/K}G\times_r H.
\end{CD}
\end{equation*}
%--------------------------------------------------------------------
We will do this in two steps,
first showing that the analogous untwisted diagram commutes, and then
showing that the twisting ideals in the various crossed products match
up properly, so that commutativity is preserved 
on taking quotients by
these ideals. 

For $N=\{e\}$, a weak form of the following theorem was proven 
in \cite[Theorem~4.1]{KQ-IC}.

\begin{thm}\label{plain-res-ind-thm}
Let $(A,G,\delta)$ be a nondegenerate
coaction, and let $N\subset H$ be closed normal
subgroups of $G$ such that Mansfield imprimitivity works for 
$H$ {\rm(}which is automatic if $H$ is amenable{\rm).}  
Then the diagram
%--------------------------------------------------------------------
\begin{equation}\label{res-ind}
\begin{CD}
\Rep A\times G/N   @>{Y_{G/N}^G}>>    \Rep A\times G\times_r N \\
@V\Res VV                       @VV\Ind V\\
\Rep A\times G/H   @>>{Y_{G/H}^G}>    \Rep A\times G\times_r H.
\end{CD}
\end{equation}
%--------------------------------------------------------------------
commutes in the usual strong sense.
\end{thm}

\begin{proof}
First note that by \cite[Theorem~5.2]{KQ-IC}, Mansfield
imprimitivity also works for $N$ and $\delta$, so $Y_{G/N}^G$ is indeed
an $A\times G\times_r N$ -- $A\times G/N$ imprimitivity bimodule, and the
above diagram makes sense. 

We shall actually prove that the Hilbert module tensor product 
\begin{equation}\label{tensor-product-eqn}
Y_{G/H}^G \otimes_{A\times G/H} \widetilde{Y_{G/N}^G} 
\end{equation}
of Mansfield bimodules is isomorphic to the reduction $X_N^H$ of
Green's bimodule, as a right-Hilbert 
$A\times G\times_r H$ -- $A\times G\times_r N$ 
bimodule.  This suffices, because then 
\begin{eqnarray*}
X_N^H \otimes_{A\times G\times_r N} Y_{G/N}^G & \cong & Y_{G/H}^G
\otimes_{A\times G/H} \widetilde{Y_{G/N}^G} \otimes_{A\times G\times_r N}
Y_{G/N}^G \\
 & \cong & Y_{G/H}^G \otimes_{A\times G/H} (A\times G/N)
\end{eqnarray*}
as a right-Hilbert $A\times G\times_r H$ -- $A\times G/N$ bimodule,
and this is exactly what it means for the above diagram to commute in
the strong sense. 

Both bimodules in the tensor product \eqref{tensor-product-eqn}\ 
are completions of Mansfield's
dense subalgebra $\D$ of $A\times G$ for the appropriate inner
products.  Our isomorphism will be the extension to 
$Y_{G/H}^G \otimes_{A\times G/H} \widetilde{Y_{G/N}^G}$ of the map
$\Phi\colon \D\otimes \tilde{\D} \to C_c(H,\D)$ defined by
$$ \Phi(x\otimes\tilde{y})(h) = x \hat{\delta}_h(y^*).$$
Note that, up to a modular function, 
$\Phi(x\otimes\tilde{y})$ is just 
Mansfield's left $C_c(H,\D)$-valued inner product
${}_{C_c(H,\D)}\langle x,y\rangle$, and hence does indeed
give an element of $C_c(H,\D)$.  
In fact, if we define $f'(h) = f(h)\Delta_H(h)^{-\frac12}$, then the
map $f\mapsto f'$ is a homeomorphism of $C_c(H,\D)$ (with the inductive
limit topology) onto itself, which takes $\Phi(\D\odot\D)$ to
${}_{C_c(H,\D)}\langle\D,\D\rangle$. This latter set is dense in
$C_c(H,A\times G)$ for the inductive limit topology
(\cite[Lemma~25]{ManJF}); it follows that the range of $\Phi$ is also
inductive limit dense in $C_c(H,A\times G)$, and therefore in $X_N^H$.
 
It only remains to show that $\Phi$ preserves the Hilbert module structure.  
For the left action of $A\times G\times_r H$, 
fix $d\in \D\subset A\times G$ and $h,t\in H$.  Then:
$$d\cdot\Phi(x\otimes\tilde{y})(h)
 =  d x \deltahat_h(y^*) 
 =  \Phi(d\cdot x\otimes\tilde{y})$$
and
\begin{eqnarray*}
t\cdot\Phi(x\otimes\tilde{y})(h) 
 & = & \deltahat_t\left(
\Phi(x\otimes\tilde{y})(t^{-1}h)\right) \Delta(t)^{\half} \\
 & = & \deltahat_t\left( x\deltahat_{t^{-1}h}(y^*)\right)
\Delta(t)^{\half} \\
 & = & \deltahat_t(x) \deltahat_h(y^*) \Delta(t)^{\half} \\
 & = & (t\cdot x) \deltahat_h(y^*) \\
 & = & \Phi( t\cdot x\otimes\tilde{y})(h).
\end{eqnarray*}
For the right action of $A\times G\times_r N$, 
fix $d\in \D\subset A\times G$, $h\in H$ and $n\in N$; 
then one similarly verifies that
\begin{eqnarray*}
\Phi(x\otimes\tilde{y})\cdot d(h) 
% & = & \Phi(x\otimes\tilde{y})(h)\deltahat_h(d) \\
% & = & x \deltahat_h(y^*d) \\
 & = & x\deltahat_h((d^*y)^*) \\
% & = & \Phi(x\otimes(d^*y)\tilde{})(h) \\
 & = & \Phi(x\otimes\tilde{y}\cdot d)(h)
\end{eqnarray*}
and
\begin{eqnarray*}
\Phi(x\otimes\tilde{y})\cdot n (h) 
% & = &
%\Phi(x\otimes\tilde{y})(hn^{-1})\, \Delta(n)^{-\half} \\
% & = & x\deltahat_{hn^{-1}}(y^*)\, \Delta(n)^{-\half} \\
 & = & x\deltahat_h(\deltahat_{n^{-1}}(y^*)\Delta(n)^{-\half}) \\
% & = & x\deltahat_h((n^{-1}\cdot y)^*) \\
 & = & \Phi(x\otimes\tilde{y}\cdot n)(h).
\end{eqnarray*}
For the right $A\times G\times_r N$-valued inner product, fix
$n\in N$ and compute:
\begin{eqnarray*}
\allowdisplaybreaks
\lefteqn{\rip{A\times G\times_r
N}{\Phi(x\otimes\tilde{y})}{\Phi(z\otimes\tilde{w})}(n) }\\
 & = & \int_H \hat{\delta}_h\left( \Phi(x\otimes\tilde{y})(h^{-1})^*
\Phi(z\otimes\tilde{w})(h^{-1}n)\right) \Delta(n)^{-{\half}}\, dh \\
 & = & \int_H \deltahat_h\left( (x\deltahat_{h^{-1}}(y^*))^*
z\deltahat_{h^{-1}n}(w^*) \right) \Delta(n)^{-\half}\, dh \\
 & = & \int_H y\deltahat_h(x^*z)\deltahat_n(w^*)\, dh\,
\Delta(n)^{-\half} \\
 & \stackrel{h\mapsto nh}{=} & \int_H y\deltahat_n\left(
\deltahat_h(x^*z)w^* \right) \, dh\, \Delta(n)^{-\half} \\
 & = & y \deltahat_n\left( w \int_H \deltahat_h(z^*x)\, dh\right)^* \,
\Delta(n)^{-\half} \\
 & = & y \deltahat_n\left( (w\cdot\rip{A\times G/H}{z}{x})^* \right) \,
\Delta(n)^{-\half} \\
 & = & \lip{A\times G\times_r N}{y}{w\cdot\rip{A\times G/H}{z}{x}}(n)
\\
 & = & \rip{A\times G\times_r
N}{x\otimes\tilde{y}}{z\otimes\tilde{w}}(n).
\end{eqnarray*}

It now follows that $\Phi$ is a 
right-Hilbert
$A\times G\times_r H$ -- $A\times G\times_r N$ bimodule isomorphism
of
$Y_{G/H}^G\otimes_{A\times G/H}\widetilde{Y_{G/N}^G}$ 
onto $X_N^H$. 
\end{proof}

\begin{cor}\label{twisted-res-ind-cor}
Let $(A,G,G/K,\delta,\tau)$ be a nondegenerate
twisted coaction, and let
$N\subset H$ be closed normal subgroups of $G$ contained in $K$ such that 
Mansfield imprimitivity works for $H$ and $\delta$
{\rm(}which is automatic if $H$ is amenable{\rm).} Then the diagram
$$
\begin{CD}
\rep A\times_{G/K}G/N @>Z^G_{G/N}>> \rep A\times_{G/K}G\times_r N \\
@V\res VV @VV\ind V \\
\rep A\times_{G/K}G/H @>>Z^G_{G/H}> \rep A\times_{G/K}G\times_r H
\end{CD}
$$
commutes in the usual strong sense.
\end{cor}

\begin{proof}
We show the appropriate ideals in diagram \eqref{res-ind} match up and
appeal to \corref{quotient}; 
because the diagram commutes, 
and the top and bottom maps are 
Morita equivalences, we need only match up the ideals
along three sides. Let $I_\tau$, $I^N_\tau$, and $I^H_\tau$ be
the twisting ideals of $A\times G$, $A\times G/N$, and $A\times G/H$,
respectively. 
The reduced crossed product $A\times_{G/K}G\times_r N$ is by definition
$\left((A\times G)/I_\tau\right)\times_r N$.  If $\pi$ is a
representation of $A\times G$ with kernel $I_\tau$ (for example, if
$\pi=k_A\times k_G$), then $\pi$ induces a faithful representation
$X_{\{e\}}^N\dashind\pi$ of $\left((A\times G)/I_\tau\right)\times_r
N$; thus
$$A\times_{G/K}\times_r N 
= (A\times G\times_r N)/(X_{\{e\}}^N\dashind I_\tau).$$
It is part of the content of \cite[Theorem~4.4]{KQ-IC}\ that 
$$Y_{G/N}^G\dashind I_\tau^N = X_{\{e\}}^N\dashind I_\tau$$
(see \cite[Equation~4.2]{KQ-IC}).  Hence the appropriate ideals match
up along the top of diagram \eqref{res-ind}, and similarly along the bottom. 

Since 
$$X^H_N\dashind(X^N_{\{e\}}\dashind I_\tau)=
X^H_{\{e\}}\dashind I_\tau,$$
by \corref{reduced-stages}, the ideals also match up along the right side of
diagram \eqref{res-ind}. 
\end{proof}

%=====================================================================

\section{Mansfield induction and Green restriction}\label{ind-res-sec}

In this section we prove analogs of \thmref{plain-res-ind-thm}\ and
\corref{twisted-res-ind-cor}, where now we use Mansfield induction on
the left sides of the diagrams and restriction on the right. 
As in the previous section, we first prove 
an untwisted version, and then show that
the twisting ideals match up properly.  

\begin{thm}\label{plain-ind-res-thm}
Let $(A,G,\delta)$ be a nondegenerate  
coaction, and let $N\subset H$ be closed normal
subgroups of $G$ such that Mansfield imprimitivity works for 
$H$ {\rm(}which is automatic if $H$ is amenable{\rm).}  
Then the diagram 
%--------------------------------------------------------------------
\begin{equation}\label{ind-res}
\begin{CD}
\Rep A\times G/N          @>Y_{G/N}^G>>    \Rep A\times G\times_r N \\
@A\Ind AA                      @AA\Res A\\
\Rep A\times G/H   @>>Y_{G/H}^G>    \Rep A\times G\times_r H.
\end{CD}
\end{equation}
%--------------------------------------------------------------------
commutes in the usual strong sense.
\end{thm}

\begin{proof}
Let us denote the Hilbert 
module for Mansfield induction from $A\times G/H$
to $A\times G/N$ by $Y_{G/H}^{G/N}$; here we are identifying $G/H$ with
$(G/N)/(H/N)$, so $Y_{G/H}^{G/N}$ is a completion of Mansfield's dense
subalgebra $\D_N$ of $A\times G/N$.  We shall prove that the Hilbert
module tensor product
$$Y_{G/N}^G \otimes_{A\times G/N} Y_{G/H}^{G/N}$$
is isomorphic to $Y_{G/H}^G$ as a right-Hilbert 
$A\times G\times_r N$ -- $A\times G/H$ 
bimodule.  This suffices, because  then
$$Y_{G/N}^G \otimes_{A\times G/N} Y_{G/H}^{G/N} \cong (A\times
G\times_r H) \otimes_{A\times G\times_r H} Y_{G/H}^G$$
as a right-Hilbert  
$A\times G\times_r N$ -- $A\times G/H$
bimodule, and this is exactly what it means for
the above diagram to commute in the strong sense. 

Our map will be 
the extension to $Y_{G/N}^G\otimes_{A\times G/N}Y_{G/H}^{G/N}$ of the
map $\Psi\colon \D\otimes \D_N\to \D$
defined by
$$\Psi(x\otimes y) = xy;$$
here the product $xy$ makes sense in $M(A\times G)$ because $y$ belongs
to $\D_N\subset M(A\times G)$. 
In other words, this product is given by the right 
action of $\D_N$ on Mansfield's bimodule $\D$,  and hence 
$\Psi(x\otimes y)$ does indeed belong to $\D$.   

We now show that $\Psi$ preserves the 
right-Hilbert bimodule structure.  Since
$A\times G$ on the left and $A\times G/H$ on the right act by
multiplication
in $M(A\times G)$, it is immediate that $\Psi$ preserves these actions. To
see that $\Psi$ preserves the left $N$-action is a straightforward
calculation, using the fact that each $y\in\D_N$ is fixed by
$\deltahat_n$ for $n\in N$:
$$\Psi(n\cdot x\otimes y)
  =  \deltahat_n(x)y\, \Delta(n)^\half 
  =  n\cdot \Psi(x\otimes y).$$
To see that $\Psi$ preserves the right $A\times G/H$-valued inner
products, note that $\widehat{\delta|}_{tN} = \deltahat_t$ on $A\times
G/N$, and compute:
\begin{eqnarray*}
\allowdisplaybreaks
\lefteqn{\rip{A\times G/H}{\Psi(x\otimes y)}{\Psi(z\otimes w)} }\\
 & = & \int_H
\deltahat_t\left( \Psi(x\otimes y)^*\Psi(z\otimes w)\right) dt \\
 & = & \int_H \deltahat_t( (xy)^*zw)\, dt \\
 & = & \int_{H/N}\int_N \deltahat_{tn}(y^*x^*zw)\, dn\, dtN \\
 & = & \int_{H/N}\int_N 
\deltahat_t( y^*\deltahat_n(x^*z) w)\, dn\, dtN \\
 & = & \int_{H/N} 
\deltahat_t \left(y^*\rip{A\times G/N}{x}{z} w\right) dtN \\
 & = & \int_{H/N}
\widehat{\delta|}_{tN} \left(y^*\rip{A\times G/N}{x}{z}\cdot w\right) dtN \\
 & = & \rip{A\times G/H}{y}{\rip{A\times G/N}{x}{z}\cdot w} \\
 & = & \rip{A\times G/H}{x\otimes y}{z\otimes w}.
\end{eqnarray*}

It only remains to show that the range of $\Psi$ is dense in
$Y_{G/H}^G$.  For this, note that $\D\cdot j_G(C_c(G/N)) = \D$, since
if $x\in \D$ is $(u,E)$, we may choose $f\in C_c(G/N)$ such that $f$ is
identically $1$ on $E$, and then $x\cdot j_G(f) = x$.  Hence, $\D\cdot
j_G(C_c(G/N))$ is dense in $Y_{G/H}^G$.  Since $Y_{G/H}^G$ is a
nondegenerate right $A$-module, and $\delta_{A_c(G)}(A)$ is dense in
$A$, we therefore have
\begin{eqnarray*}
\overline{\Phi(\D\odot\D_N)} & = & \overline{\D\cdot \D_N}\\
 & = & \overline{\D\cdot j_G(C_c(G/N))j_A(\delta_{A_c(G)}(A))}\\
 & = & \overline{\D\cdot j_A(\delta_{A_C(G)}(A))}\\
 & = & Y_{G/H}^G.
\end{eqnarray*}
It now follows that $\Psi$ is a 
right-Hilbert
$A\times G\times_r N$ -- $A\times G/H$ bimodule isomorphism
of
$Y_{G/N}^G\otimes_{A\times G/N}Y_{G/H}^{G/N}$ onto 
$Y_{G/H}^G$. 
\end{proof}

\begin{cor}\label{mansfield-stages-thm}
{\rm(}Mansfield Induction in Stages{\rm.)}
Let $(A,G,\delta)$ be a nondegenerate
coaction, and let $N\subset H$ be closed normal
subgroups of $G$ such that Mansfield imprimitivity works for
$H$ {\rm(}which is automatic if $H$ is amenable{\rm).}  Then 
the diagram 
\begin{equation*}
%----------------------------------------------------------------------
\begin{diagram}
	\node{\rep A\times G/N}
		\arrow{e,t}{Y_{G/N}^G}
	\node{\rep A\times G}\\
	\node{\rep A\times G/H}
		\arrow{n,l}{Y_{G/H}^{G/N}}
		\arrow{ne,b}{Y_{G/H}^G}
\end{diagram}
%----------------------------------------------------------------------
\end{equation*}
commutes in the usual strong sense.
\end{cor}

\begin{proof}
The corollary requires that
the Hilbert
module tensor product
$$Y_{G/N}^G \otimes_{A\times G/N} Y_{G/H}^{G/N}$$
be isomorphic to $Y_{G/H}^G$ as right-Hilbert 
$A\times G$ -- $A\times G/H$ bimodules; we 
showed slightly more than this in 
the proof of \thmref{plain-ind-res-thm}\ by including
the left $N$-action.
\end{proof}

\begin{cor}\label{twisted-ind-res-cor}
Let $(A,G,G/K,\delta,\tau)$ be a nondegenerate
twisted coaction, 
and let $N\subset H$ be closed
normal subgroups of $G$ contained in $K$ such that
Mansfield imprimitivity works for $H$ and $\delta$
{\rm(}which is automatic if $H$ is amenable{\rm)}.
Then the diagram
$$
\begin{CD}
\rep A\times_{G/K}G/N @>Z^G_{G/N}>> \rep A\times_{G/K}G\times_r N \\
@A\ind AA @AA\res A \\
\rep A\times_{G/K}G/H @>>Z^G_{G/H}> \rep A\times_{G/K}G\times_r H
\end{CD}
$$
commutes in the usual strong sense.
\end{cor}

\begin{proof}
As with \corref{twisted-res-ind-cor}, we need only show the appropriate
ideals in the diagram \eqref{ind-res} match up, and appeal to
\corref{quotient}. As before, 
let $I_\tau$, $I^N_\tau$, and $I^H_\tau$ be
the twisting ideals of $A\times G$, $A\times G/N$, and $A\times G/H$,
respectively.
The reduced crossed product $A\times_{G/K}G\times_r N$ is by definition
$\left((A\times G)/I_\tau\right)\times_r N$.  If $\pi$ is a
representation of $A\times G$ with kernel $I_\tau$ (for example, if
$\pi=k_A\times k_G$), then $\pi$ induces a faithful representation
$X_{\{e\}}^N\dashind\pi$ of $\left((A\times G)/I_\tau\right)\times_r
N$; thus
$$A\times_{G/K}\times_r N
= (A\times G\times_r N)/(X_{\{e\}}^N\dashind I_\tau).$$
It is part of the content of \cite[Theorem~4.4]{KQ-IC}\ that
$$Y_{G/N}^G\dashind I_\tau^N = X_{\{e\}}^N\dashind I_\tau$$
(see \cite[Equation~4.2]{KQ-IC}).  Hence the appropriate ideals match
up along the top of diagram \eqref{ind-res}, and similarly along the
bottom.

It only remains to see that the ideals match up along the right side
of \eqref{ind-res}.
We need
\begin{equation}
\res_{A\times G\times_r H}^{A\times G\times_r N}
\ind^{A\times G\times_r H}_{A\times G}I_\tau=
\ind^{A\times G\times_r N}_{A\times G}I_\tau.
\end{equation}
Now, our reduced version \corref{reduced-stages}\ 
of induction in stages gives
\begin{equation}
\ind^{A\times G\times_r H}_{A\times G}I_\tau=
\ind^{A\times G\times_r H}_{A\times G\times_r N}
\ind^{A\times G\times_r N}_{A\times G}I_\tau,
\end{equation}
so it will suffice to show two things:
\begin{enumerate}
\item the ideal $\ind^{A\times G\times_r N}_{A\times G}
I_\tau$ of $A\times G\times_r
N$ is $H$-invariant;
\item if $J$ is an $H$-invariant ideal of $A\times G\times_r N$, then
$$\res_{A\times G\times_r H}^{A\times G\times_r N}
\ind^{A\times G\times_r H}_{A\times G\times_r N}J=J.$$
\end{enumerate}
These are shown abstractly in the next two lemmas, which
although at most partly new, may have
independent interest.
\end{proof}

\begin{lem}
Let $(B,H,\alpha)$ be an action, $N$ a closed normal subgroup of $H$, and $I$
an $H$-invariant ideal of $B$. Then the ideal 
$\ind_B^{B\times_r N}I$ of
$B\times_r N$ is $H$-invariant for the decomposition action.
\end{lem}

\begin{proof}
The decomposition action of $H$ on $B\times N$ leaves the kernel of the
regular representation invariant, hence it indeed induces an action
$\beta$ of $H$ on $B\times_r N$.  Explicitly, for $n\in N$, $s\in H$,
and $c\in C_c(N,B)$, 
$$\beta_s(c)(n)=\gamma(s)\alpha_s(c(s^{-1}ns)),$$
where $\gamma$ is the modular function of conjugation of $H$ on
$N$:
$$\int_N f(n)\,dn=\gamma(s)\int_N f(s^{-1}ns)\,dn
\righttext{for}f\in C_c(N),s\in H.$$

Let $X=X^N_{\{e\}}$ be the usual right-Hilbert 
$B\times_r N$ -- $B$ 
bimodule, so $X$ is a completion of $C_c(N,B)$. We will show $X$ is
$H$-equivariant, which will imply that $H$-invariant ideals of $B$ induce to
$H$-invariant ideals of $B\times_r N$. More precisely, we will construct
a strongly continuous Banach representation $u$ of $H$ on $X$ such that
for
$s\in H$,
$x,y\in X$, $b\in B$, and $c\in B\times_r N$ we have
\begin{gather}
\langle u_s(x),u_s(y)\rangle_B \label{isometric}
=\alpha_s(\langle x,y\rangle_B)\\
u_s(xb)=u_s(x)\alpha_s(b) \label{rightaction}\\
u_s(cx)=\beta_s(c)u_s(x). \label{leftaction}
\end{gather}
A straightforward calculation then shows that if $I$ is $H$-invariant,
so is $X\dashind I = \{ c\in B\times_r N\mid \rip{B}{c\cdot
x}{y}\in I \text{ for all } x,y\in X\}$.  

For $s\in H$, $x\in C_c(N,B)$ define $u_s(x)\in C_c(N,B)$ by
$$u_s(x)(n)=\gamma(s)^{\half}\alpha_s(x(s^{-1}ns))
\righttext{for}n\in N.$$
To show \eqref{isometric}, take $x,y\in C_c(N,B)$ and compute:
\begin{align}
\begin{split}
\langle u_s(x),u_s(y)\rangle_B
&=\int_N \alpha_n\bigl(u_s(x)(n^{-1})^*u_s(y)(n^{-1})\bigr)\,dn \\
&=\int_N \gamma(s)
\alpha_{ns}\bigl(x(s^{-1}n^{-1}s)^*y(s^{-1}n^{-1}s)\bigr)\,dn \\
&=\int_N \alpha_{sn}\bigl(x(n^{-1})^*y(n^{-1}\bigr)\,dn \\
&=\alpha_s\left(\int_N \alpha_n\bigl(
x(n^{-1})^*y(n^{-1}\bigr)\,dn\right) \\
&=\alpha_s(\langle x,y\rangle_B).
\end{split}
\end{align}
It is now clear that $u$ is a homomorphism of $H$ into the isometric
automorphisms of the normed space $C_c(N,B)$, hence determines by
continuity a homomorphism of $H$ into the isometric automorphisms of the
Banach space $X$, and \eqref{isometric} follows, again by continuity.

For \eqref{rightaction}, it suffices to take $x\in C_c(N,B)$: for $n\in
N$ we have
\begin{align}
\begin{split}
u_s(xb)(n)
&=\gamma(s)^{\half}\alpha_s\bigl((xb)(s^{-1}ns)\bigr)\\
&=\gamma(s)^{\half}\alpha_s\bigl(x(s^{-1}ns)\alpha_{s^{-1}ns}(b)\bigr)\\
&=\gamma(s)^{\half}\alpha_s(x(s^{-1}ns))\alpha_{ns}(b)\\
&=u_s(x)(n)\alpha_n(\alpha_s(b))\\
&=\bigl(u_s(x)\alpha_s(b)\bigr)(n).
\end{split}
\end{align}

We next show $u$ is strongly continuous.
By \cite[\S2]{GreLS}, it suffices to show that $u$ is strongly
continuous for the inductive limit topology on $C_c(N,B)$.  Fix $x\in
C_c(N,B)$ with compact support $E$, and suppose $s_i\to e$ in $G$.
Since $u_{s_i}(x)\to x$ uniformly (by a standard compactness argument),
we need only find a compact set $F$ in $N$ and $k$ such that
the support of $u_{s_i}(x) - x$ is in $F$ for $i\geq k$.  To do this,
choose a neighborhood $U$ of $e$ in $G$ with compact closure, and let
$k$ be such that $s_i\in U$ for $i\geq k$.  Then the compact set $F = N\cap
(\bar{U}E\bar{U}^{-1})$ will do.  

Finally, for \eqref{leftaction}, it suffices to take $c,x\in C_c(N,B)$:
for $n\in N$ we have
\begin{align}
\begin{split}
u_s(bx)(n)
&=\gamma(s)^{\half}\alpha_s\bigl((cx)(s^{-1}ns)\bigr)\\
&=\gamma(s)^{\half}\alpha_s\left(
\int_N c(k)\alpha_k(x(k^{-1}s^{-1}ns))\,dk\right)\\
&=\gamma(s)^{\half}\int_N \alpha_s(c(k))
\alpha_{sk}(x(k^{-1}s^{-1}ns))\,dk\\
&=\gamma(s)^{\half}\int_N \gamma(s)\alpha_s(c(s^{-1}ks))
\alpha_{ks}(x(s^{-1}k^{-1}ns))\,dk\\
&=\int_N \beta_s(c)(k)\alpha_k(u_s(x)(k^{-1}n))\,dk\\
&=\bigl(\beta_s(c)u_s(x)\bigr)(n).
\end{split}
\end{align}
\end{proof}

\begin{lem}\label{res-ind-invt-ideal-lem}
If $(B,H,\alpha)$ is an action, $N$ is a 
closed normal subgroup of $H$, and $J$
is an $H$-invariant ideal of $B\times_r N$, then
$$\res_{B\times_r H}^{B\times_r N}
\ind^{B\times_r H}_{B\times_r N}J=J.$$
\end{lem}

\begin{proof}
We first show the analogous equality
\begin{equation}\label{invariant}
\res_{B\times H}^{B\times N}
\ind^{B\times H}_{B\times N}K=K
\end{equation}
for full crossed products, assuming $K$ is an $H$-invariant ideal of
$B\times N$. By \cite[Proposition 1]{GreLS}, the decomposition action
of $H$ on $B\times N$ is twisted over $N$, and
$$B\times H\cong B\times N\times_N H.$$
So \eqref{invariant} follows from \cite[Proposition 11]{GreLS}.

Now we have an $H$-equivariant commutative diagram
$$
\begin{CD}
B\times N @>>> M(B\times H) \\
@V{\rho_N}VV @VV{\rho_H}V \\
B\times_r N @>>> M(B\times_r H),
\end{CD}
$$
where $\rho_N$ denotes the regular representation of $B\times N$. 
This gives an $H$-equivariant commutative diagram
$$
\begin{CD}
\rep B\times N @<{\res_{B\times H}^{B\times N}}<< \rep B\times H \\
@A{\rho^*_N}AA @AA{\rho^*_H}A \\
\rep B\times_r N @<{\res_{B\times_r H}^{B\times_r N}}<< \rep B\times_r H.
\end{CD}
$$
Also, we have a diagram
$$
\begin{CD}
\rep B\times N @>{\ind^{B\times H}_{B\times N}}>> \rep B\times H \\
@A{\rho^*_N}AA @AA{\rho^*_H}A \\
\rep B\times_r N @>{\ind^{B\times_r H}_{B\times_r N}}>> \rep B\times_r H,
\end{CD}
$$
which commutes by our construction of the bottom arrow in 
\lemref{reduced}. Since
$J$ is an
$H$-invariant ideal of $B\times_r N$, $\rho^*_N J$ is an $H$-invariant ideal
of $B\times N$, and we have
\begin{align}
\begin{split}
\rho^*_N \res_{B\times_r H}^{B\times_r N}
\ind^{B\times_r H}_{B\times_r N}J
&=\res_{B\times H}^{B\times N} \rho^*_H
\ind^{B\times_r H}_{B\times_r N}J \\
&=\res_{B\times H}^{B\times N}
\ind^{B\times H}_{B\times N} \rho^*_N J \\
&=\rho^*_N J
\end{split}
\end{align}
by \eqref{invariant}. Since $\rho_N\colon B\times N\to B\times_r N$ is
surjective, $\rho^*_N$ is injective on ideals, and the lemma follows.
\end{proof}

%=====================================================================

\section{Morita equivalence, inflation, and stabilization}
\label{me-infl-stab-sec}

In this section we show that our Res-Ind duality is compatible with
certain standard constructions. 

Before discussing Morita equivalence of coactions, we recall 
the concept of multiplier bimodules introduced in
\cite{ER-MI}.
A \emph{multiplier} $m=(m_A,m_B)$ of ${}_AX_B$
consists of an $A$-linear map $m_A\colon A\to X$ and a $B$-linear map
$m_B\colon B\to X$ such that 
$m_A(a)b=am_B(b)$ for $a\in A$, $b\in B$. The \emph{multiplier
bimodule} $M(X)$ consists of all multipliers of $X$. 
An \emph{imprimitivity bimodule homomorphism}
$\phi=(\phi_A,\phi_X,\phi_B)\colon {}_AX_B\to M({}_CY_D)$ consists of
homomorphisms $\phi_A\colon A\to M(C)$ and $\phi_B\colon B\to M(D)$ and
a bimodule map $\phi_X\colon X\to M(Y)$ preserving the inner products:
\begin{gather*}
{}_{M(C)}\langle\phi_X(x),\phi_X(y)\rangle
=\phi_A({}_A\langle x,y\rangle)
\quad\text{and} \\
\langle\phi_X(x),\phi_X(y)\rangle_{M(D)}
=\phi_B(\langle x,y\rangle_B)
\quad\text{for }x,y\in X.
\end{gather*}
$\phi$ is called \emph{nondegenerate} if $\phi_A$ and $\phi_B$ are
nondegenerate.
Significantly, this implies an ostensibly stronger form of nondegeneracy:

\begin{lem}\label{iains-ib-hom-lem}
If $\phi\colon {}_AX_B\to M({}_CY_D)$ is a nondegenerate imprimitivity
bimodule homomorphism, then
\[
\overline{C\phi_X(X)}=Y=\overline{\phi_X(X)D}.
\]
\end{lem}

\begin{proof}
By symmetry, it suffices to show the second equality. Since
$Y=\overline{C\cdot Y}$, and $\phi_A$ is nondegenerate, we have
$Y=\overline{\phi_A(A)\cdot Y}$. Now, since the range of the inner
product ${}_A\langle\cdot,\cdot\rangle$ spans $A$, we have
\[
\begin{split}
Y
&=\overline{\phi_A({}_A\langle X,X\rangle)\cdot Y} \\
&=\overline{{}_{M(C)}\langle \phi_X(X),\phi_X(X)\rangle \cdot Y} \\
&=\overline{\phi_X(X)\cdot\langle \phi_X(X),Y\rangle_{M(D)}}. 
\end{split}
\]
Because we can factor $Y=\overline{C\cdot Y}$, and $C\cdot M(Y)\subset
Y$, the pairing $\langle\cdot,\cdot\rangle_{M(D)}$ takes $M(Y)\times Y$
to $D$, and we can deduce from our calculation that $Y\subset
\overline{\phi_X(X)\cdot D}$. The other inclusion is trivial, so this
establishes the lemma.
\end{proof}

\begin{lem}
\label{lem:nondegenerate}
Let ${}_AX_B$ and ${}_CY_D$ be imprimitivity bimodules, and let
$\phi\colon L(X)\to M(L(Y))$ be a nondegenerate homomorphism such
that
\[
\phi\begin{pmatrix}1 & 0 \\ 0 & 0\end{pmatrix}
=\begin{pmatrix}1 & 0 \\ 0 & 0\end{pmatrix},
\]
where $L(X)=\left(\begin{smallmatrix}A & X
\\ \tilde X & B\end{smallmatrix}\right)$ is the linking algebra of $X$,
and similarly for $Y$. Then $\phi$ restricts on the corners to give a
nondegenerate imprimitivity bimodule homomorphism $\Phi\colon {}_AX_B\to
M({}_CY_D)$.
\end{lem}

\begin{proof}
Since $\phi$ is nondegenerate, we have $\phi\left(\begin{smallmatrix}0 &
0 \\ 0 & 1\end{smallmatrix}\right)=\left(\begin{smallmatrix}0 & 0 \\ 0 &
1\end{smallmatrix}\right)$. Since
\[
M(L(Y))=\left(\begin{smallmatrix}M(C)
& M(Y) \\ M(\tilde Y) & M(D)\end{smallmatrix}\right),
\]
we deduce that
there are unique linear maps
\begin{gather*}
\Phi_A\colon A\to M(C),\quad \Phi_X\colon X\to M(Y), \\
\Phi_{\tilde X}\colon \tilde X\to M(\tilde Y),\quad\text{and}\quad
\Phi_B\colon B\to M(D)
\end{gather*}
such that
\[
\phi=\begin{pmatrix}\Phi_A & \Phi_X \\ \Phi_{\tilde X} &
\Phi_B\end{pmatrix}.
\]
The algebraic properties of $\phi$ show that
$\Phi:=(\Phi_A,\Phi_X,\Phi_B)$ is an imprimitivity bimodule
homomorphism of ${}_AX_B$ to $M({}_CY_D)$, and nondegeneracy of $\Phi_A$
and $\Phi_B$ is inherited from nondegeneracy of $\phi$.
\end{proof}

We will need the following elementary fact about imprimitivity bimodule
homomorphisms.

\begin{lem}\label{surjection}
If $\phi = (\phi_A,\phi_X,\phi_B)\colon {}_AX_B\to M({}_CY_D)$ is a
nondegenerate imprimitivity bimodule homomorphism, then the diagram
$$
\begin{CD}
\rep A @<X<< \rep B \\
@A{\phi^*_A}AA @AA{\phi^*_B}A \\
\rep C @<<Y< \rep D
\end{CD}
$$
commutes in the usual strong sense.
\end{lem}

\begin{proof}
We must show $X\otimes_BD\cong C\otimes_C Y$ as a right-Hilbert
$A$ -- $D$ bimodule. 
Of course, ${}_AC\otimes_C Y$ is isomorphic to ${}_AY$.
Define $\Psi\colon
X\odot D\to Y$ by
$$\Psi(x\otimes d)=\phi_X(x)d.$$
The required properties of $\Psi$ follow from the nondegeneracy of
$\phi$, and straightforward calculations showing that
$$\langle \Psi(x\otimes d),\Psi(y\otimes e)\rangle_D
= \langle x\otimes d,y\otimes e\rangle_D$$
and
$$\Psi(a(x\otimes d)e)
=a\,\Psi(x\otimes d)e$$
for $a\in A$, $x,y\in X$, and $d,e\in D$. 
\end{proof}

Following
\cite[Definition~3.3]{NgCC}
(see also \cite{BS-CH}, \cite{BuiTC}, \cite{ER-MI}),
a \emph{coaction}\ $\delta$ of $G$ on an
imprimitivity bimodule ${}_AX_B$ is an imprimitivity bimodule
homomorphism
$$\delta =
(\delta_A,\delta_X,\delta_B)\colon {}_AX_B\to M\left({}_{A\otimes
C^*(G)}(X\otimes C^*(G))_{B\otimes C^*(G)}\right)$$
such that $(A,G,\delta_A)$ and $(B,G,\delta_B)$ are $C^*$-coactions,
and
satisfying
$$(\delta_X\otimes\id)\circ\delta_X =
(\id\otimes\delta_G)\circ\delta_X.$$
As a consequence of the definition, we automatically have
$\delta_X(x)\cdot(1_B\otimes z)$ and $(1_X\otimes z)\cdot\delta_X(x)
\in X\otimes C^*(G)$ for $x\in X$, $z\in C^*(G)$.  Also, since by
assumption $\delta_A$ and $\delta_B$ are nondegenerate
$C^*$-homomorphisms, $\delta$ is automatically  nondegenerate
as an imprimitivity bimodule homomorphism.  Hence, by
\lemref{iains-ib-hom-lem}, we have
$$\overline{\delta_X(X)\cdot(A\otimes C^*(G))} 
= \overline{(B\otimes C^*(G))\cdot\delta_X(X)} = X\otimes C^*(G).$$
When such a $\delta$ exists we say $(A,G,\delta_A)$ and
$(B,G,\delta_B)$ are \emph{Morita equivalent}, and we call
$(X,\delta_X)$
a \emph{Morita equivalence} of $\delta_A$ and $\delta_B$.

If $N$ is a closed normal subgroup of $G$, then
$$\delta_X|=(\iota\otimes q_N)\circ\delta_X\colon 
X\to M(X\otimes C^*(G/N))$$
is a coaction of $G/N$ on $X$,
where $q_N\colon C^*(G)\to C^*(G/N)$ is the canonical
quotient map. A \emph{Morita equivalence} between
twisted coactions $(A,G,G/K,\delta_A,\tau_A)$ and
$(B,G,G/K,\delta_B,\tau_B)$ is an $(A,G,\delta_A)$ -- $(B,G,\delta_B)$
Morita equivalence $(X,\delta_X)$ such that
$$\delta_X|(x)=\tau_A\otimes\iota(w_{G/K})(x\otimes 1)
\tau_B\otimes\iota(w_{G/K})^*\righttext{for}x\in X.$$
In this case, for any closed normal subgroup $N$ of $G$ 
contained in $K$, there
are $A\times_{G/K}G$ -- $B\times_{G/K}G$ and
$A\times_{G/K}G/N$ -- $B\times_{G/K}G/N$ imprimitivity bimodules
$X\times_{G/K}G$ and $X\times_{G/K}G/N$, respectively. Moreover, there
is an action of $K$ on $X\times_{G/K}G$ inducing a Morita equivalence
between the dual actions $(A\times_{G/K}G,K,\tilde\delta_A)$ 
and $(B\times_{G/K}G,K,\tilde\delta_B)$, so
there is an $A\times_{G/K}G\times_r K$ -- $B\times_{G/K}G\times_r K$
imprimitivity bimodule $X\times_{G/K}G\times_r K$.

The next result shows \corref{twisted-res-ind-cor}\ is 
compatible with Morita equivalence:

\begin{thm}\label{me-res-ind-cube-thm}
If $(X,\delta_X)$ is a Morita equivalence between nondegenerate twisted
coactions $(A,G,G/K,\delta_A,\tau_A)$ and $(B,G,G/K,\delta_B,\tau_B)$,
and $N\subset H$ are closed normal subgroups of $G$ contained in $K$
such that Mansfield
imprimitivity works for $H$ and one of the
coactions, then the cube
%--------------------------------------------------------------------
\begin{equation}\label{res-ind ME}
\tiny
\dgARROWLENGTH=0em
\begin{diagram}
\node[2]{\Rep B\times_{G/K}G/N}
        \arrow[2]{e}
        \arrow{s,-}
	\arrow{sw,t}{X\times_{G/K}G/N}
\node[2]{\Rep B\times_{G/K}G\times_r N}
        \arrow[2]{s,r}{\Ind}
        \arrow{sw,t}{X\times_{G/K}G\times_r N}\\
\node{\Rep A\times_{G/K}G/N}
        \arrow[2]{e}
        \arrow[2]{s,l}{\Res}
\node{}
        \arrow{s,l}{\Res}
\node{\Rep A\times_{G/K}G\times_r N}
        \arrow[2]{s,r,3}{\Ind}\\
\node[2]{\Rep B\times_{G/K}G/H}
        \arrow{e,-}
        \arrow{sw,b}{X\times_{G/K}G/H}
\node{}
        \arrow{e}
\node{\Rep B\times_{G/K}G\times_r H}
        \arrow{sw,b}{X\times_{G/K}G\times_r H}\\
\node{\Rep A\times_{G/K}G/H}
        \arrow[2]{e}
\node[2]{\Rep A\times_{G/K}G\times_r H}
\end{diagram}
\end{equation}
%--------------------------------------------------------------------
commutes in the usual strong sense.
\end{thm}

\begin{proof}
First note that by \cite[Theorem~5.3]{KQ-IC}, if Mansfield
imprimitivity works for $H$ and one of the coactions, it works for the
other, so the above cube makes sense.

Because all the horizontal arrows are bijections, 
we need only show commutativity of three of the vertical faces,
as well as the top and bottom. The front and back
commute by 
\corref{twisted-res-ind-cor}.

The bottom face is just the top face with $N$ replaced by $H$; 
we show the top commutes. The
untwisted version (i.e., with $K=G$)
%---------------------------------------------------------------------
\begin{equation}\label{untwisted res-ind ME}
\begin{diagram}
\node{\Rep B\times_{}G/N}
        \arrow{e}
\node{\Rep B\times_{}G\times_r N}\\
\node{\Rep A\times_{}G/N}
        \arrow{n,l}{X\times_{}G/N}
        \arrow{e}
\node{\Rep A\times_{}G\times_r N}
        \arrow{n,r}{X\times_{}G\times_r N}
\end{diagram}
\end{equation}
%---------------------------------------------------------------------
is \cite[Proposition 4.5]{ER-ST}. Even though they use reduced
coactions and require $N$ to be amenable, their arguments carry over to
our setting since we assume that Mansfield imprimitivity works,
as pointed out in \cite{KQ-IC}. 
So, it only remains to show the ideals in diagram 
\eqref{untwisted res-ind ME}\ match up and appeal to 
\corref{quotient}.
The ideals match up along the top and bottom 
by \cite[Theorem~4.4]{KQ-IC}, and along the left side by
\cite[Corollary 3.3]{ER-MI}, and this is enough, since the top and
bottom are Morita equivalences.

For the left face, we use the following general lemma:

\begin{lem}\label{co-res-res-lem}
Let $(X,\delta_X)$ be a Morita equivalence between coactions
$(A,G,\delta_A)$ and $(B,G,\delta_B)$, and let $N$ be a closed normal
subgroup of $G$. Then the diagram
%--------------------------------------------------------------------
\begin{equation}\label{co-res-res-diag}
\begin{diagram}
\node{\rep A\times G}
	\arrow{s,l}{\res}
\node{\rep B\times G}
	\arrow{w,t}{X\times G}
	\arrow{s,r}{\res} \\
\node{\rep A\times G/N}
\node{\rep B\times G/N}
	\arrow{w,b}{X\times G/N}
\end{diagram}
\end{equation}
%--------------------------------------------------------------------
commutes.
\end{lem}

\begin{proof}
We aim to apply \lemref{surjection}, so we need a nondegenerate imprimitivity
bimodule homomorphism
\[
\Phi\colon {}_{A\times G/N}X\times G/N_{B\times
G/N}\to M({}_{A\times G}X\times G_{B\times G})
\]
with
\[
\Phi_{A\times G/N}=j_A\times j^A_{G}|
\quad\text{and}\quad
\Phi_{B\times G/N}=j_B\times j^B_{G}|.
\]
Let $L(X)$ be the linking algebra for $X$; then there
is a coaction $\left(\begin{smallmatrix}\delta_A & \delta_X \\
\delta_{\tilde X} & \delta_B\end{smallmatrix}\right)$
of $G$ on $L(X)$. By [ER, Appendix], we
have
\[
L(X)\times G/N=L(X\times G/N)
\]
and
\[
M(L(X)\times G))=\begin{pmatrix}M(A\times G) & M(X\times G) \\
M(\tilde X\times G) & M(B\times G)\end{pmatrix}.
\]
By \lemref{lem:nondegenerate}, the nondegenerate homomorphism
$$j_{L(X)}\times j^{L(X)}_G|\colon L(X)\times G/N\to M(L(X)\times G)$$
restricts on the corners to a nondegenerate imprimitivity bimodule
homomorphism 
\[
\Phi\colon {}_{A\times G/N}X\times G/N_{B\times G/N}\to
M({}_{A\times G}X\times G_{B\times G}).
\]
 Since the restrictions of
$j_{L(X)}\times j^{L(X)}_G|$ to the diagonal corners $A\times G/N$ and
$B\times G/N$ agree with $j_A\times j^A_G|$ and $j_B\times j^B_G|$,
respectively, we are done.
\end{proof}

Returning to the proof of 
\thmref{me-res-ind-cube-thm}, the untwisted version of the left-hand
face of diagram
\eqref{res-ind ME}\ follows from \lemref{co-res-res-lem}.  
We
need only show that the 
twisting ideals in 
each
crossed product in diagram \eqref{co-res-res-diag}\ match up, and apply
\corref{quotient}.  
That the ideals match up across the horizontal arrows is shown in
\cite[Corollary~3.3]{ER-MI}.  That the ideals match up along the
vertical arrows follows from the proof of \corref{twisted-res-ind-cor}.
This completes the proof of \thmref{me-res-ind-cube-thm}
\end{proof}

In the proof of \thmref{me-res-ind-cube-thm}, we were able to deduce
commutativity of the right face of \eqref{res-ind ME} from the
other five faces.
This is a special case of a general ``$\Ind$-$\Ind$'' diagram 
for actions which is related to \cite[Theorem~3]{EchME}\ and 
\cite[Proposition~4.1.2]{KalME}.

We now show that \corref{twisted-ind-res-cor}\ is 
compatible with Morita equivalence:

\begin{thm}\label{me-ind-res-cube-thm}
If $(X,\delta_X)$ is a Morita equivalence between nondegenerate twisted
coactions $(A,G,G/K,\delta_A,\tau_A)$ and $(B,G,G/K,\delta_B,\tau_B)$,
and $N\subset H$ are closed normal subgroups of $G$ contained in $K$
such that Mansfield
imprimitivity works for $H$ and one of the
coactions, then the cube
%---------------------------------------------------------------------
\begin{equation}\label{ind-res ME}
\tiny
\dgARROWLENGTH=0em
\begin{diagram}
\node[2]{\Rep B\times_{G/K}G/N}
        \arrow[2]{e}
        \arrow{sw,t}{X\times_{G/K}G/N}
\node[2]{\Rep B\times_{G/K}G\times_r N}
        \arrow{sw,t}{X\times_{G/K}G\times_r N}\\
\node{\Rep A\times_{G/K}G/N}
        \arrow[2]{e}
\node{}
        \arrow{n,l}{\Ind}
\node{\Rep A\times_{G/K}G\times_r N}\\
\node[2]{\Rep B\times_{G/K}G/H}
        \arrow{n,-}
        \arrow{e,-}
	\arrow{sw,b}{X\times_{G/K}G/H}
\node{}
        \arrow{e}
\node{\Rep B\times_{G/K}G\times_r H}
        \arrow[2]{n,r}{\Res}
        \arrow{sw,b}{X\times_{G/K}G\times_r H}\\
\node{\Rep A\times_{G/K}G/H}
        \arrow[2]{n,l}{\Ind}
        \arrow[2]{e}
\node[2]{\Rep A\times_{G/K}G\times_r H}
        \arrow[2]{n,r,3}{\Res}
\end{diagram}
\end{equation}
%---------------------------------------------------------------------
commutes in the usual strong sense.
\end{thm}

\begin{proof}
The front and back faces commute by 
\corref{twisted-ind-res-cor}, and the top and bottom faces are
the same as in \eqref{res-ind ME}, so it suffices to show commutativity
of the right face.
For this we use the following lemma:

\begin{lem}\label{rt-face-res-lem}
Let $(Y,\gamma)$ be a Morita equivalence between actions $(C,H,\alpha)$
and $(D,H,\beta)$, and let $N$ be a closed subgroup of $H$.
Then the diagram
%---------------------------------------------------------------------
\begin{equation}\label{act-res-res-diag}
\begin{diagram}
\node{\rep C\times H}
\arrow{s,l}{\res}
\node{\rep D\times H}
\arrow{w,t}{Y\times H}
\arrow{s,r}{\res} \\
\node{\rep C\times N}
\node{\rep D\times N}
\arrow{w,b}{Y\times N}
\end{diagram}
\end{equation}
%---------------------------------------------------------------------
commutes.
\end{lem}

\begin{proof}
We aim to apply \lemref{surjection}, 
so we need a nondegenerate imprimitivity
bimodule homomorphism
\[
\Phi\colon {}_{C\times N}Y\times N_{D\times N}\to
M({}_{C\times H}Y\times H_{D\times H})
\]
with
\[
\Phi_{C\times N}=i_C\times i^C_H|N
\quad\text{and}\quad
\Phi_{D\times N}=i_D\times i^D_H|N.
\]
Let $L(Y)$ be the linking algebra for $Y$; then there
is an action $\left(\begin{smallmatrix}\alpha & \gamma \\ \tilde\gamma &
\beta\end{smallmatrix}\right)$ of $H$ on $L(Y)$. By [ER, Appendix], we
have
\[
L(Y)\times N=L(Y\times N)
\]
and
\[
M(L(Y)\times H))=\begin{pmatrix}M(C\times H) & M(Y\times H) \\
M(\tilde Y\times H) & M(D\times H)\end{pmatrix}.
\]
By \lemref{lem:nondegenerate}, the nondegenerate homomorphism
\[
i_{L(Y)}\times i^{L(Y)}_H|N\colon L(Y)\times N\to M(L(Y)\times H)
\]
restricts on the corners to a nondegenerate imprimitivity bimodule
homomorphism 
\[
\Phi\colon {}_{C\times N}Y\times N_{D\times N}\to
M({}_{C\times H}Y\times H_{D\times H}).
\]
Since the restrictions of
$i_{L(Y)}\times i^{L(Y)}_H|N$ to the diagonal corners $C\times N$ and
$D\times N$ agree with $i_C\times i^C_H|N$ and $i_D\times i^D_H|N$,
respectively, we are done.
\end{proof}

Returning to the proof of the right-hand face of diagram
\eqref{ind-res ME}, we
need only show that the kernels of the regular representations of each
crossed product in diagram \eqref{act-res-res-diag}\ match up, and apply
\corref{quotient}.  
That the kernels match up across the horizontal
arrows is shown in \cite[\S6]{ComCP}. 
The kernels match up along the vertical
arrows because the regular representation of
$N$ is quasi-equivalent to the regular representation of $H$
restricted to $N$.  Alternatively, we can argue as follows: we must
show that
$$ \Res^{C\times N}_{C\times H} \Ind_C^{C\times H} \{0\} =
\Ind_C^{C\times N} \{0\}.$$
By Green induction in stages,
$$ \Ind_C^{C\times H}  = \Ind_{C\times N}^{C\times H}\circ
\Ind_C^{C\times N}.$$
Since for $N$ normal, $\Ind_C^{C\times N}\{0\}$ is $H$-invariant,
equation \eqref{invariant}\ gives us
\begin{eqnarray*}
\Res^{C\times N}_{C\times H} \Ind_C^{C\times H} \{0\}
 & = & \Res^{C\times N}_{C\times H} \Ind_{C\times N}^{C\times H}
\Ind_C^{C\times N}  \{0\}\\
 & = & \Ind_C^{C\times N}  \{0\}.
\end{eqnarray*}
This completes the proof of \thmref{me-ind-res-cube-thm}.
\end{proof}

\begin{rem}
By \cite[Proposition~2.3]{KQ-IC}, 
Morita equivalence respects nondegeneracy of coactions.  
Hence, in Theorems~\ref{me-res-ind-cube-thm}\
and \ref{me-ind-res-cube-thm}, one coaction is in fact nondegenerate if
and only if the other one is.  
\end{rem}

This time, we were able to
deduce commutativity of the left face of \eqref{ind-res ME} from
the other five faces. Written in terms of tensor products,
\begin{multline}
Z^{G/N}_{G/H}(A)\otimes_{A\times_{G/K}G/H}
(X\times_{G/K}GH)\\
\cong(X\times_{G/K}G/N)\otimes_{B\times_{G/K}G/N}
Z^{G/N}_{G/H}(B),
\end{multline}
where we have used self-explanatory notation to distinguish the 
bimodules for $A$ and $B$. The special case $N=\{e\}$ is \cite[Theorem
4.4]{ER-ST}. A direct proof is presumably possible,
although probably quite tedious.

We turn to inflation of coactions: if $K$ is a closed normal subgroup of $G$
and $(A,K,\epsilon)$ is a coaction, composing with
the natural embedding of $A\otimes C^*(K)$
in $M(A\otimes C^*(G))$ gives a coaction $(A,G,
\inf\epsilon)$, called an \emph{inflated} coaction. $\inf\epsilon$ is
trivially twisted over $G/K$ by $f\mapsto f(e)1$, and
\cite[Example 2.14]{PR-CO}\ gives a natural isomorphism of 
$A\times_{\inf\epsilon,G/K}G$ onto $A\times_{\epsilon}K$ which takes
$k_A(a)k_G(f)$ to $j_A(a)j_K(f|_K)$.  

The next two theorems show Corollaries \ref{twisted-res-ind-cor}\ and 
\ref{twisted-ind-res-cor}\ are compatible with
inflation:

\begin{thm}\label{inf-res-ind-cube-thm}
If $N\subset H$ are closed normal subgroups of $G$ contained in $K$, and
$(A,K,\epsilon)$ is a nondegenerate
coaction such that Mansfield imprimitivity works for 
$H$ and $\epsilon$, 
then the cube
%----------------------------------------------------------------------
\begin{equation}\label{res-ind inf}
\tiny
\dgARROWLENGTH=0em
\begin{diagram}
\node[2]{\Rep A\times K/N}
        \arrow[2]{e}
        \arrow{s,-}
\node[2]{\Rep A\times K\times_r N}
        \arrow[2]{s,r}{\Ind}\\
\node{\Rep A\times_{G/K}G/N}
        \arrow{ne,t}{\cong}
        \arrow[2]{e}
        \arrow[2]{s,l}{\Res}
\node{}
        \arrow{s,l}{\Res}
\node{\Rep A\times_{G/K}G\times_r N}
        \arrow{ne,t}{\cong}
        \arrow[2]{s,r,3}{\Ind}\\
\node[2]{\Rep A\times K/H}
        \arrow{e,-}
\node{}
        \arrow{e}
\node{\Rep A\times K\times_r H}\\
\node{\Rep A\times_{G/K}G/H}
        \arrow{ne,b}{\cong}
        \arrow[2]{e}
\node[2]{\Rep A\times_{G/K}G\times_r H}
        \arrow{ne,b}{\cong}
\end{diagram}
\end{equation}
%----------------------------------------------------------------------
commutes in the usual strong sense.
\end{thm}

\begin{proof}
First, by \cite[Theorem~5.4]{KQ-IC}, Mansfield imprimitivity
works for $H$ and the 
(nondegenerate) inflated coaction $(A,G,\inf\epsilon)$, so the
above diagram makes sense. 

The front and back faces commute by \thmref{plain-res-ind-thm}.

We show the top commutes; the bottom is essentially the
same. \cite[Proposition 4.8
and Remark 4.9]{ER-ST}\ give an imprimitivity bimodule homomorphism
$$(\Phi_N,\Phi_Y,\Phi_{G/N})\colon Y^G_{G/N}\to Y^K_{K/N}.$$
Hence, \lemref{surjection}\ gives a commutative diagram
%---------------------------------------------------------------------
\begin{equation}\label{inflated}
\begin{CD}
\rep A\times G/N @>{Y^G_{G/N}}>> \rep A\times G\times_r N \\
@A{\Phi_{G/N}^*}AA @AA{\Phi_N^*}A \\
\rep A\times K/N @>>{Y^K_{K/N}}> \rep A\times K\times_r N.
\end{CD}
\end{equation}
%---------------------------------------------------------------------
To pass from \eqref{inflated} to the top of \eqref{res-ind inf}, we need
to show the appropriate ideals $I^N_\tau$ in $A\times G/N$,
$X^N_{\{e\}}\dashind I_\tau$ in $A\times G\times_r N$, and $\{0\}$ in
both $A\times K/N$ and $A\times K\times_r N$ match up. The ideals match
up along the top (by \cite[Theorem~4.4]{KQ-IC}), 
bottom (trivially), and left
(since $\Phi_{G/N}(I^N_\tau)=\{0\}$), and this is enough.

Finally, for the left face of \eqref{res-ind inf}, the diagram
%---------------------------------------------------------------------
\[
\begin{diagram}
\node{A\times_{G/K}G/N}
\arrow{e,t}{\cong}
\node{A\times K/N} \\
\node{A\times_{G/K}G/H}
\arrow{n}
\arrow{e,t}{\cong}
\node{A\times K/H}
\arrow{n}
\end{diagram}
\]
%---------------------------------------------------------------------
of nondegenerate homomorphisms commutes, since the diagram
%---------------------------------------------------------------------
\[
\begin{diagram}
\node{M(C_0(G/N))}
\arrow{e}
\node{M(C_0(K/N))} \\
\node{C_0(G/H)}
\arrow{n}
\arrow{e}
\node{C_0(K/H)}
\arrow{n}
\end{diagram}
\]
%---------------------------------------------------------------------
does. Hence, the required diagram
%---------------------------------------------------------------------
\[
\begin{diagram}
\node{\rep A\times_{G/K}G/N}
\arrow{e,t}{\cong}
\arrow{s,l}{\res}
\node{\rep A\times K/N}
\arrow{s,r}{\res} \\
\node{\rep A\times_{G/K}G/H}
\arrow{e,b}{\cong}
\node{\rep A\times K/H}
\end{diagram}
\]
%---------------------------------------------------------------------
of right-Hilbert bimodules commutes as well.
\end{proof}

\begin{thm}\label{inf-ind-res-cube-thm}
If $N\subset H$ are closed normal subgroups of $G$ contained in $K$, and
$(A,K,\epsilon)$ is a nondegenerate
coaction such that Mansfield imprimitivity works for 
$H$ and $\epsilon$, then the cube
%---------------------------------------------------------------------
\begin{equation}\label{ind-res inf}
\tiny
\dgARROWLENGTH=0em
\begin{diagram}
\node[2]{\Rep A\times K/N}
        \arrow[2]{e}
\node[2]{\Rep A\times K\times_r N}\\
\node{\Rep A\times_{G/K}G/N}
        \arrow{ne,t}{\cong}
        \arrow[2]{e}
\node{}
        \arrow{n,l}{\Ind}
\node{\Rep A\times_{G/K}G\times_r N}
        \arrow{ne,t}{\cong}\\
\node[2]{\Rep A\times K/H}
        \arrow{n,-}
        \arrow{e,-}
\node{}
        \arrow{e}
\node{\Rep A\times K\times_r H}
        \arrow[2]{n,r}{\Res}\\
\node{\Rep A\times_{G/K}G/H}
        \arrow[2]{n,l}{\Ind}
        \arrow{ne,b}{\cong}
        \arrow[2]{e}
\node[2]{\Rep A\times_{G/K}G\times_r H}
        \arrow[2]{n,r,3}{\Res}
        \arrow{ne,b}{\cong}
\end{diagram}
\end{equation}
%---------------------------------------------------------------------
commutes in the usual strong sense.
\end{thm}

\begin{proof}
The front and back faces commute by 
\thmref{plain-ind-res-thm}, and the top and bottom faces are
the same as in \eqref{res-ind inf}.  Since the 
horizontal maps in the
right face come from equivariant isomorphisms, the commutativity of
this face follows from \lemref{rt-face-res-lem}\ as in the proof of
\thmref{me-ind-res-cube-thm}. 
\end{proof}

\begin{rem}
By \cite[Proposition~2.1]{KQ-IC}, in Theorems~\ref{inf-res-ind-cube-thm}\
and \ref{inf-ind-res-cube-thm}, $\inf\epsilon$ is in fact nondegenerate
if and only if $\epsilon$ is.  Moreover, when the coactions are
nondegenerate, Mansfield imprimitivity works for $H$ and $\inf\epsilon$
if and only if it works for $H$ and $\epsilon$.  
\end{rem}

As before, we get commutativity of a sixth face of diagrams 
\eqref{res-ind inf} and \eqref{ind-res inf} from the other five faces.
The left face 
of \eqref{ind-res inf} is new; the special case $N=\{e\}$ follows from
\cite[Theorem 4.7]{ER-ST}.

Finally, we show that our Res-Ind duality is compatible with the
stabilization trick of \cite{ER-ST}, as adapted to full coactions in
\cite{KQ-IC}.  
Let $(A,G,G/K,\delta,\tau)$ be a nondegenerate twisted
coaction such that Mansfield imprimitivity works for 
$\delta$ and $K$ itself;
equivalently, such that $\delta$ is 
normal \cite[Lemma~3.6]{KQ-IC}.
Then Mansfield imprimitivity works for $\delta$ and any closed normal
subgroup of $G$ contained in $K$ by \cite[Theorem~5.2]{KQ-IC}\ (or by
\cite[Lemma~3.2]{KQ-IC}, since $\delta$ is normal). 
\cite[Theorem 3.1]{ER-ST}\
shows that 
the twisted coaction $(A,G,G/K,\delta,\tau)$ is Morita equivalent to  
the (nondegenerate) inflated twisted coaction $(A\times_{G/K}G\times_r
K,G,G/K,\inf(\widehat{\widetilde\delta})^{\rm n},1)$, so by \cite[Theorems 5.3
and 5.4]{KQ-IC}, Mansfield imprimitivity works for
$(\widehat{\widetilde\delta})^{\rm n}$ and 
any closed normal subgroup of $G$
contained in $K$. 
Thus we can chain together 
Theorems \ref{me-res-ind-cube-thm}\ and \ref{inf-res-ind-cube-thm},
and similarly Theorems \ref{me-ind-res-cube-thm}\ and
\ref{inf-ind-res-cube-thm}, to obtain:

\begin{thm}\label{stab-cube-thm}
Let $(A,G,G/K,\delta,\tau)$ be a nondegenerate normal twisted coaction, 
and let $N\subset H$ be closed normal subgroups of $G$ contained in $K$. 
Further let $B=A\times_{G/K}G\times_rK$, which carries the double dual
coaction of $K$. 
Then the cubes
%---------------------------------------------------------------------
\begin{equation}\label{res-ind stab}
\tiny
\dgARROWLENGTH=0em
\begin{diagram}
\node[2]{\Rep B \times K/N}
	\arrow[2]{e}
	\arrow{s,-}
\node[2]{\Rep B\times K\times_r N}
	\arrow[2]{s,r}{\Ind}\\
\node{\Rep A\times_{G/K}G/N}
	\arrow{ne}
	\arrow[2]{e}
	\arrow[2]{s,l}{\Res}
\node{}
	\arrow{s,l}{\Res}
\node{\Rep A\times_{G/K}G\times_r N}
	\arrow{ne}
	\arrow[2]{s,r,3}{\Ind}\\
\node[2]{\Rep B\times K/H}
	\arrow{e,-}
\node{}
	\arrow{e}
\node{\Rep B\times K\times_r H}\\
\node{\Rep A\times_{G/K}G/H}
	\arrow{ne}
	\arrow[2]{e}
\node[2]{\Rep A\times_{G/K}G\times_r H}
	\arrow{ne}
\end{diagram}
\end{equation}
%---------------------------------------------------------------------
and
%---------------------------------------------------------------------
\begin{equation}\label{ind-res stab}
\tiny
\dgARROWLENGTH=0em
\begin{diagram}
\node[2]{\Rep B\times K/N}
        \arrow[2]{e}
\node[2]{\Rep B\times K\times_r N}\\
\node{\Rep A\times_{G/K}G/N}
        \arrow{ne}
        \arrow[2]{e}
\node{}
        \arrow{n,l}{\Ind}
\node{\Rep A\times_{G/K}G\times_r N}
        \arrow{ne}\\
\node[2]{\Rep B\times K/H}
	\arrow{n,-}
        \arrow{e,-}
\node{}
        \arrow{e}
\node{\Rep B\times K\times_r H}
	\arrow[2]{n,r}{\Res}\\
\node{\Rep A\times_{G/K}G/H}
	\arrow[2]{n,l}{\Ind}
        \arrow{ne}
        \arrow[2]{e}
\node[2]{\Rep A\times_{G/K}G\times_r H}
	\arrow[2]{n,r,3}{\Res}
        \arrow{ne}
\end{diagram}
\end{equation}
%---------------------------------------------------------------------
both commute in the usual strong sense.
\end{thm}

Again, the left face of \eqref{ind-res stab}\ is new; when $N=\{e\}$ it
reduces to \cite[Theorem 4.7]{ER-ST}.

%====================================================================

\section{Ind, Res, Ex, Sub}\label{ind-res-ex-sub-sec}

In this section, as a sample application of our $\Res$-$\Ind$ duality,
we generalize some results of Gootman and Lazar \cite[\S3]{GL-AN}\
concerning restriction and induction of ideals in crossed products by
coactions of amenable groups, to nonamenable groups.  
Nilsen \cite{NilRE}\ has recently proved similar results, using
different, representation-theoretic techniques.
Our methods, based on our Res-Ind duality results, appear to be more
efficient than those of Gootman and Lazar.
With some additional effort, we could further 
generalize to the setting of intermediate
twisted crossed products by coactions of nonamenable groups, but we
feel that to do so at this point would only muddy the waters. 

First, recall from \cite[Proposition~9]{GreLS}\ that a nondegenerate
homomorphism $\phi\colon A\to M(B)$ gives rise to maps 
$\Res_\phi=\phi^*\colon \I(B)\to \I(A)$ 
and
$\Ex_\phi=\phi_*\colon \I(A)\to \I(B)$ 
between spaces of ideals.  
By definition, 
$\Res_\phi(J) = \{ a\in A \mid \phi(a)B\subset J\}$, 
and $\Ex_\phi(I)$ is the ideal generated by $\phi(I)B$.  
If in addition $X$ is a $B$ -- $C$ imprimitivity bimodule,
then by definition we have the following commutative diagram:
%----------------------------------------------------------------------
\begin{equation}
\begin{diagram}
	\node{\I(B)}
		\arrow{s,l}{\Res_\phi}
	\node{\I(C)}
		\arrow{w,t}{X\dashind}
		\arrow{sw,b}{\Ind_\phi}\\
	\node{\I(A).}
\end{diagram}
\end{equation}
%----------------------------------------------------------------------
Green defines a map $\Sub_\phi$ by requiring the diagram
%----------------------------------------------------------------------
\begin{equation}
\begin{diagram}
	\node{\I(B)}
		\arrow{e,t}{\tilde{X}\dashind}
	\node{\I(C)}\\
	\node{\I(A)}
		\arrow{n,l}{\Ex_\phi}
		\arrow{ne,b}{\Sub_\phi}
\end{diagram}
\end{equation}
%----------------------------------------------------------------------
to commute.  
(Actually, we have abstracted Green's definition a bit in order to 
bring the properties of the maps into high relief.) 
So, $\Sub$ is to $\Ind$ as $\Ex$ is to $\Res$, and in fact
$\Ex$ is a special case of $\Sub$, just as $\Res$ is a special case of
$\Ind$.  

We define the \emph{sup} of a set of ideals to be the ideal they
generate (i.e., the closed span of the union of the ideals). 
In the above situation we have the following facts: 

\begin{prop}\cite[Proposition~9(i)]{GreLS}\label{GAN-prop}
Fix ideals $I$ of $A$, $J$ of $B$, and $K$ of $C$.  Then:
\begin{itemize}
\openup1\jot
\item[(i)] $\Res\Ex I\contains I$; $\Ind\Sub I\contains I${\rm;}
\item[(ii)] $\Ex\Res J\subset J$; $\Sub\Ind K\subset K${\rm;}
\item[(iii)] $\Ind$, $\Res$, $\Sub$ and $\Ex$ are order-preserving
\item[(iv)] $\Ind$ and $\Res$ preserve arbitrary intersections{\rm;}
\item[(v)] $\Sub$ and $\Ex$ preserve arbitrary sups{\rm;}
\item[(vi)] $\Ex\Res\Ex I = \Ex I$; $\Sub\Ind\Sub I = \Sub I${\rm;}
\item[(vii)] $\Res\Ex\Res J = \Res J$; $\Ind\Sub\Ind K = \Ind K${\rm;}
\item[(viii)] $\Ex I$ is the smallest ideal $J'$ of $B$ such that $\Res
J'\contains I${\rm;}
\item[(ix)] $\Res J$ is the largest ideal $I'$ of $A$ such that $\Ex
I'\subset J${\rm;}
\item[(x)] $\Sub I$ is the smallest ideal $K'$ of $C$ such that $\Ind
K'\contains I${\rm;}
\item[(xi)] $\Ind K$ is the largest ideal $I'$ of $A$ such that $\Sub
I'\subset K.$
\end{itemize}
\end{prop}

Now consider an action $(B,G,\alpha)$.  Green applies the above
machinery to the canonical map $i_B\colon B\to M(B\times G)$ to get
maps
$$\Res = \Res_{i_B}\colon \I(B\times G)\to \I(B)$$
and
$$\Ex = \Ex_{i_B}\colon \I(B)\to \I(B\times G).$$
He then applies the same machinery to the canonical map $j_{B\times
G}\colon B\times G\to M(B\times G\times G)$ and the $B\times G\times G$
-- $B$ imprimitivity bimodule $X$ from his imprimitivity theorem to get
maps
$$\Ind = 
\Res_{j_{B\times G}}\circ X\dashind\colon \I(B)\to \I(B\times G)$$
and
$$\Sub = 
\tilde{X}\dashind\circ\Ex_{j_{B\times G}}\colon \I(B\times G)\to
\I(B).$$
(Green uses $C_0(G,B)\times G$ as the imprimitivity algebra, rather
than the isomorphic co-crossed product $B\times G\times G$.)

An ideal $J\subset B$ is \emph{$G$-invariant} if
$\alpha_s(J)\subset J$ for each $s\in G$;
we denote the $G$-invariant ideals of $B$ by ${}^G\I(B)$.  Green proves:

\begin{prop}\cite[Proposition~11]{GreLS}\label{green-sum-prop}
Let $(B,G,\alpha)$ be an action, and let $J$ be an ideal of $B$.  Then:
\begin{itemize}
\openup1\jot
\item[(i)] $\Res J$ and $\Sub J$ are $G$-invariant{\rm;}
\item[(ii)] $\Res\Ind J$ is the largest $G$-invariant ideal of $B$
contained in $J${\rm;}
\item[(iii)] $\Res\Ex J$ is the smallest $G$-invariant ideal of
$B$ containing $J${\rm;}
\item[(iv)] $J$ is $G$-invariant if and only if $\Res\Ind J = J$.
\end{itemize}
\end{prop}

We will have to adapt Green's machinery to
reduced crossed products.  Let $(B,G,\alpha)$ be an
action, and let $\rho\colon B\times G\to B\times_{r}G$ be the regular
representation.  The commutative diagram
%----------------------------------------------------------------------
\begin{equation}
\begin{diagram}
	\node{B}
		\arrow{e,t}{i_B}
		\arrow{se,b}{\rho\circ i_B}
	\node{M(B\times G)}
		\arrow{s,r}{\rho}\\
	\node[2]{M(B\times_{r}G)}
\end{diagram}
\end{equation}
%----------------------------------------------------------------------
gives in turn commutative diagrams
%----------------------------------------------------------------------
\begin{equation}
\begin{diagram}
	\node{\I(B)}
		\arrow{e,t}{\Ex_{i_B}}
		\arrow{se,b}{\Ex_{\rho\circ i_B}}
	\node{\I(B\times G)}
		\arrow{s,r}{\rho_*}\\
	\node[2]{\I(B\times_{r}G)}
\end{diagram}
\end{equation}
%----------------------------------------------------------------------
and
%----------------------------------------------------------------------
\begin{equation}\label{res-res-q-diag}
\begin{diagram}
	\node{\I(B)}
	\node{\I(B\times G)}
		\arrow{w,t}{\Res_{i_B}}\\
	\node[2]{\I(B\times_{r}G),}
		\arrow{nw,b}{\Res_{\rho\circ i_B}}
		\arrow{n,r}{\rho^*}
\end{diagram}
\end{equation}
%----------------------------------------------------------------------
The horizontal 
$\Ex$ and $\Res$ are what Green uses; we will need the diagonal
ones, together with their   
associated maps $\Sub_{\rho\circ i_B}\colon \I(B\times_r
G)\to \I(B)$ and $\Ind_{\rho\circ i_B}\colon \I(B)\to \I(B\times_r G)$. 

Similarly, the commutative diagram
%---------------------------------------------------------------------
\begin{equation}
\begin{diagram}
	\node{B\times G}
		\arrow{s,l}{\rho}
		\arrow{e,t}{j_{B\times G}}
	\node{M(B\times G\times G)}
		\arrow{s,r}{\cong}\\
	\node{B\times_{r}G}
		\arrow{e,b}{j_{B\times_r G}}
	\node{M(B\times_{r}G\times G)}
\end{diagram}
\end{equation}
%---------------------------------------------------------------------
gives in turn commutative diagrams
%---------------------------------------------------------------------
\begin{equation}\label{sub-5-diag}
\begin{diagram}
	\node{\I(B\times G\times G)}
		\arrow[2]{s,l}{\cong}
		\arrow{se,b}{\cong}
	\node[2]{\I(B\times G)}
		\arrow[2]{w,t}{\Ex}
		\arrow{sw,b}{\Sub}
		\arrow[2]{s,r}{\rho_*}\\
	\node[2]{\I(B)}
		\arrow{sw,t}{\cong}\\
	\node{\I(B\times_{r}G\times G)}
	\node[2]{\I(B\times_{r}G)}
		\arrow[2]{w,b}{\Ex}
		\arrow{nw,t}{\Sub}
\end{diagram}
\end{equation}
%---------------------------------------------------------------------
and
%---------------------------------------------------------------------
\begin{equation}\label{ind-5-diag}
\begin{diagram}
	\node{\I(B\times G\times G)}
		\arrow[2]{s,l}{\cong}
		\arrow{se,b}{\cong}
		\arrow[2]{e,t}{\Res}
	\node[2]{\I(B\times G)}\\
	\node[2]{\I(B)}
		\arrow{sw,t}{\cong}
		\arrow{ne,b}{\Ind}
		\arrow{se,t}{\Ind}\\
	\node{\I(B\times_{r}G\times G)}
		\arrow[2]{e,b}{\Res}
	\node[2]{\I(B\times_{r}G).}
		\arrow[2]{n,r}{\rho^*}
\end{diagram}
\end{equation}
%---------------------------------------------------------------------

Diagrams \eqref{res-res-q-diag}\ and \eqref{ind-5-diag}, together with the
$G$-equivariance and injectivity of 
$\rho^*$, imply that all of Green's results
in \cite[Lemma~10(ii) and Proposition~11]{GreLS}\
carry over to reduced crossed products.  In particular, those
results summarized in \propref{green-sum-prop}\ carry over; we will
cite these as simply \cite{GreLS}\ without further comment. 

We now consider a 
nondegenerate normal 
coaction $(A,G,\delta)$.  We would like to apply the
general abstract nonsense of the beginning of this section to obtain
maps $\Sub$ and $\Ex$ among the spaces $\I(A)$ and
$\I(A\times G)$  to go along with  the maps $\Res$ and
$\Ind$. We use the canonical map 
$j_A\colon A\to M(A\times G)$ 
to get a map
$$\Ex = \Ex_{j_A}\colon \I(A)\to \I(A\times G).$$
(In fact, we've already used this map in diagram
\eqref{sub-5-diag}\ in the case $A=B\times G$.) 
Similarly, we use the canonical map 
$i_{A\times G}\colon A\times G\to M(A\times G\times_r G)$ 
and the $A\times G\times_r G$
-- $A$ imprimitivity bimodule $Y$ from Mansfield's
imprimitivity theorem to get a map
$$\Sub =
\tilde{Y}\dashind\circ\Ex_{i_{A\times G}}\colon \I(A\times G)\to
\I(A).$$

In this situation, \thmref{plain-res-ind-thm}\ in the case $N=\{e\}$
and $H=G$ gives us the commutativity of the following diagram:
%---------------------------------------------------------------------
\begin{equation}\label{res-ind-triangle}
\begin{diagram}
	\node{\I(A\times G)}
		\arrow{s,l}{\Res}
		\arrow{se,t}{\Ind}\\
	\node{\I(A)}
		\arrow{e,b}{Y\dashind}
	\node{\I(A\times G\times_r G).}
\end{diagram}
\end{equation}
%---------------------------------------------------------------------
By the same token, \thmref{plain-ind-res-thm}\ reduces to the
commutativity of
%---------------------------------------------------------------------
\begin{equation}\label{ind-res-triangle}
\begin{diagram}
	\node{\I(A\times G)}\\
	\node{\I(A)}
		\arrow{n,l}{\Ind}
		\arrow{e,b}{Y\dashind}
	\node{\I(A\times G\times_r G),}
		\arrow{nw,t}{\Res}
\end{diagram}
\end{equation}
%---------------------------------------------------------------------
which is just the definition of $\Ind$.  

We can deduce dualities between $\Sub$ and $\Ex$
from Corollaries 
\ref{twisted-res-ind-cor}\ and \ref{twisted-ind-res-cor}.

\begin{thm}\label{sub-ex-thm}
Let $(A,G,G/K,\delta,\tau)$ be a nondegenerate twisted coaction, and
let $N\subset H$ be closed normal subgroups of $G$ contained in $K$ such that
Mansfield imprimitivity works for $H$ and $\delta$.  
Then the following diagrams commute{\rm:}
%---------------------------------------------------------------------
\begin{equation}\label{twisted-sub-ex-diag}
\begin{diagram}
	\node{\I(A\times_{G/K}G/N)}
		\arrow{s,l}{\Sub}
		\arrow{e,t}{Z_{G/N}^G\dashind}
	\node{\I(A\times_{G/K}G\times_r N)}
		\arrow{s,r}{\Ex}\\
	\node{\I(A\times_{G/K}G/H)}
		\arrow{e,b}{Z_{G/H}^G\dashind}
	\node{\I(A\times_{G/K}G\times_r H)}
\end{diagram}
\end{equation}
%---------------------------------------------------------------------
%---------------------------------------------------------------------
\begin{equation}\label{twisted-ex-sub-diag}
\begin{diagram}
	\node{\I(A\times_{G/K}G/N)}
		\arrow{e,t}{Z_{G/N}^G\dashind}
	\node{\I(A\times_{G/K}G\times_r N)}\\
	\node{\I(A\times_{G/K}G/H)}
		\arrow{n,l}{\Ex}
		\arrow{e,b}{Z_{G/H}^G\dashind}
	\node{\I(A\times_{G/K}G\times_r H).}
		\arrow{n,r}{\Sub}
\end{diagram}
\end{equation}
%---------------------------------------------------------------------
\end{thm}

\begin{proof}
The theorem follows immediately from Corollaries 
\ref{twisted-res-ind-cor}\ and \ref{twisted-ind-res-cor}, 
together with elementary properties of $\Res$, $\Ex$, $\Sub$ and
$\Ind$.  For example, by \cite[\S3]{RieUR},
\corref{twisted-ind-res-cor}\ gives a commutative diagram
%---------------------------------------------------------------------
\begin{equation}
\begin{diagram}
        \node{\I(A\times_{G/K}G/N)}
                \arrow{e,t}{Z_{G/N}^G\dashind}
        \node{\I(A\times_{G/K}G\times_r N)}\\
        \node{\I(A\times_{G/K}G/H)}
                \arrow{n,l}{\Ind}
                \arrow{e,b}{Z_{G/H}^G\dashind}
        \node{\I(A\times_{G/K}G\times_r H).}
                \arrow{n,r}{\Res}
\end{diagram}
\end{equation}
%---------------------------------------------------------------------
Now fix $I\in \I(A\times_{G/K}G/N)$.  By \propref{GAN-prop}, $\Sub I$
is the smallest ideal $J$ of $A\times_{G/K}G/H$ such that $\Ind
J\contains I$; hence $Z_{G/H}^G\dashind(\Sub I)$ is the smallest ideal
$K$ of $A\times_{G/K}G\times_r H$ such that $\Res K\contains
Z_{G/N}^G\dashind I$.  Again using \propref{GAN-prop}, this implies
that 
$$Z_{G/H}^G\dashind(\Sub I) = \Ex(Z_{G/N}^G\dashind I),$$ 
so diagram
\eqref{twisted-sub-ex-diag}\ commutes.  

A similar argument works for diagram~\eqref{twisted-ex-sub-diag}. 
\end{proof}

For the rest of this section, we will only need the special case
$N=\{e\}$ and $H=K=G$ of \thmref{sub-ex-thm}.  For reference, this
gives the commutativity of the following triangles:
%---------------------------------------------------------------------
\begin{equation}\label{sub-ex-triangle}
\begin{diagram}
        \node{\I(A\times G)}
                \arrow{s,l}{\Sub}
                \arrow{se,t}{\Ex}\\
        \node{\I(A)}
                \arrow{e,b}{Y\dashind}
        \node{\I(A\times G\times_r G)}
\end{diagram}
\end{equation}
%---------------------------------------------------------------------
%---------------------------------------------------------------------
\begin{equation}\label{ex-sub-triangle}
\begin{diagram}
        \node{\I(A\times G)}\\
        \node{\I(A)}
                \arrow{n,l}{\Ex}
                \arrow{e,b}{Y\dashind}
        \node{\I(A\times G\times_r G).}
                \arrow{nw,t}{\Sub}
\end{diagram}
\end{equation}
%---------------------------------------------------------------------
(Again, diagram \eqref{sub-ex-triangle}\ has simply reduced to the
definition of $\Sub$.) 

We are now able to generalize some of the coaction results of
\cite[\S3]{GL-AN}.  Let $(A,G,\delta)$ be a nondegenerate normal coaction. 
Adapting \cite[Definition~2.4]{ER-ST}\ to full coactions,
we say that an ideal $I$ of $A$ is
\emph{$G$-invariant} if 
$$I = \ker (q_I\otimes\lambda)\circ\delta,$$
where $q_I\colon A\to A/I$ is the quotient map
and $\lambda\colon C^*(G)\to C^*_r(G)$ is the left regular
representation.  We do not know if this
implies that $\delta$ restricts to a coaction on $I$ (unless $G$ is
amenable), but 
we do get a coaction $\delta_{A/I}$ of $G$ on
$A/I$.  Since we will not need this fact here, we omit the proof.

Let ${}^G\I(A)$ denote the $G$-invariant ideals of $A$. 

\begin{lem}\label{res-ind-invt-lem}
Let $(A,G,\delta)$ be a nondegenerate normal coaction. 
\begin{itemize}
\openup1\jot
\item[(i)] An ideal $I$ of $A$ is $G$-invariant if and only if 
$$I = \Res\Ind I.$$
\item[(ii)] $\Ind$ is injective when 
restricted to ${}^G\I(A)$, and $\Res$ is onto ${}^G\I(A)$. 
\end{itemize}
\end{lem}

\begin{proof}
Part (ii) is immediate from (i).  To see (i), we have:
\begin{eqnarray*}
\ker (q_I\otimes\lambda)\circ\delta & = & \ker( \Res\Ind q_I)\\
 & = & \Res\Ind(\ker q_I) \\
 & = & \Res\Ind I.
\end{eqnarray*}
\end{proof}

Green shows \cite[Proposition~13]{GreLS}\ that for an action
$(B,G,\alpha)$ with $G$ amenable, $\Ex = \Ind$ on ${}^G\I(B)$.  
This is definitely not true for non-amenable $G$: $\Ex \{0\}
= \{0\}$, while $\Ind\{0\} = \{0\}$ if and only if $B\times G =
B\times_{r}G$.  
For all we know, even if we pass to reduced crossed products,
$\Ex$ and $\Ind$ can be different on ${}^G\I(B)$.
Even for trivial actions, the question reduces to the unsolved problem
of whether there exists a locally compact group $G$ for which $C^*_r(G)$
is not exact.

Coactions behave like actions of abelian groups, so
the following generalization of
\cite[Proposition~3.14(iii)]{GL-AN}\ is not surprising:

\begin{lem}\label{ex=ind-lem} 
Let $(A,G,\delta)$ be a 
nondegenerate normal coaction, and let $I\in {}^G\I(A)$.  Then
$$\Ex I = \Ind I.$$
\end{lem}

\begin{proof}
Using the dualities \eqref{ind-res-triangle}, \eqref{res-ind-triangle},
and \eqref{ex-sub-triangle}, together with \propref{green-sum-prop}\
(i) and (iv), we have:
\begin{eqnarray*}
\Ind\Res\Ex I & = & \Res\Ind\Sub Y\dashind I\\
 & = & \Sub Y\dashind I\\
 & = & \Ex I,
\end{eqnarray*}
so since $\Ind$ is order-preserving and $\Res\Ex$ is increasing, we have 
$$\Ind I\subset \Ind\Res\Ex I = \Ex I.$$

On the other hand, by invariance we have $I = \Res\Ind I$, hence
$$\Ex I = \Ex\Res\Ind I \subset \Ind I,$$
since $\Ex\Res$ is decreasing.  
\end{proof}

We do not know how to prove
the companion result for $\Sub$ and $\Res$
(generalizing \cite[Proposition~3.14(iv)]{GL-AN}); it certainly doesn't
follow from \cite{GreLS}\ and duality as in the
proof of \lemref{ex=ind-lem}.

We next show that the maps $\Res$, $\Ind$ and $\Ex$
produce invariant ideals, extending parts of 
\cite[Propositions 3.14(i) and
3.15(i)]{GL-AN}\ to the non-amenable case.  

\begin{prop}\label{invt-ideals-prop}
Let $(A,G,\delta)$ be a nondegenerate normal
coaction, and fix $I\in \I(A)$ and
$J\in\I(A\times G)$.  Then{\rm:}
\begin{itemize}
\openup1\jot
\item[(i)] $\Ind I$ and $\Ex I$ are in ${}^G\I(A\times G)${\rm;}
\item[(ii)] $\Res J$ is in ${}^G\I(A)$.  
\end{itemize}
\end{prop}

\begin{proof}
By the $\Res$-$\Ind$ and $\Sub$-$\Ex$ dualities
\eqref{res-ind-triangle}\ and \eqref{sub-ex-triangle}, 
$$\Ind I = \Res Y\dashind I$$
and 
$$\Ex I = \Sub Y\dashind I,$$
so (i) follows from \cite[Proposition~11(i)]{GreLS}. 

For (ii), let ${}_GJ$ denote the largest $G$-invariant ideal 
contained in $J$.  Then \cite[Proposition~11(ii)]{GreLS}\ gives
\begin{eqnarray*}
\Res\Ind\Res J & = & \tilde{Y}\dashind \Ind\Res\Ind J\\
 & = & \tilde{Y}\dashind \Ind {}_GJ\\
 & = & \tilde{Y}\dashind \Ind J\\
 & = & \Res J, 
\end{eqnarray*}
so $\Res J$ is invariant by \lemref{res-ind-invt-lem}. 
\end{proof}

The following proposition generalizes
\cite[Proposition~3.14(v)]{GL-AN}.

\begin{prop}\label{iri-rir-prop}
Let $(A,G,\delta)$ be a nondegenerate normal
coaction, and fix $I\in \I(A)$ and
$J\in\I(A\times G)$.  Then $\Ind\Res\Ind I = \Ind I$ 
and $\Res\Ind\Res J = \Res J$.  
\end{prop}

\begin{proof}
Since $\Ind I$ is invariant, we have
$$\Ind I = \Res\Ind\Ind I = \Ind\Res\Ind I$$
by duality.

Since $\Res J$ is invariant, we immediately have
$$\Res\Ind\Res J = \Res J.$$
\end{proof}

The following lemma generalizes 
part of \cite[Lemma~3.12]{GL-AN}, which relies
explicitly on the amenability of $G$. 

\begin{lem}\label{closure-lem}
Let $(A,G,\delta)$ be a nondegenerate normal
coaction.  Then ${}^G\I(A)$ is closed under
arbitrary intersections. 
\end{lem}

\begin{proof}
$\Res$ is onto ${}^G\I(A)$ (\lemref{res-ind-invt-lem}) and preserves
arbitrary intersections; the lemma follows. 
\end{proof}

Our final result generalizes part of
\cite[Proposition~3.15, (iii) and
(iv)]{GL-AN}:

\begin{prop}\label{max-min-prop}
Let $(A,G,\delta)$ be a nondegenerate normal
coaction, and fix $I\in \I(A)$ and
$J\in\I(A\times G)$.  Then{\rm:}
\begin{itemize}
\openup1\jot
\item[(i)] $\Res\Ex I$ is the smallest $G$-invariant ideal of $A$
containing $I${\rm;}
\item[(ii)] $\Ind\Res J$ is the largest $G$-invariant ideal of
$A\times G$ contained in $J${\rm;}
\item[(iii)] $\Ind\Sub J$ is the smallest $G$-invariant ideal of
$A\times G$ containing $J$.
\end{itemize}
\end{prop}

\begin{proof}
(i) Let $K$ be the intersection of all invariant ideals of $A$
containing $I$; since ${}^G\I(A)$ is closed under intersections, $K$ is
invariant and hence is the smallest such ideal of $A$.  Now $\Res\Ex I$
is invariant, and $\Res\Ex I\contains I$ because $\Res\Ex$ is
increasing; so we get
$$I\subset K\subset \Res\Ex I.$$
Taking $\Ex$ and using \propref{GAN-prop}(vi), we have
$$\Ex I \subset \Ex K\subset \Ex\Res\Ex K = \Ex K,$$
so equality holds throughout.  In particular, 
$$\Ex K = \Ex\Res\Ex I;$$
since $\Ex$ is injective on invariant ideals, we get
$$\Res\Ex I = K.$$

(ii) By dualities \eqref{ind-res-triangle}\ and
\eqref{res-ind-triangle}, 
$$\Ind\Res J = \Res\Ind J,$$
which is the largest $G$-invariant ideal of $A\times G$ contained in
$J$, by \cite[Proposition~11(ii)]{GreLS}.

(iii) By dualities \eqref{ind-res-triangle}\ and
\eqref{sub-ex-triangle}, 
$$\Ind\Sub J = \Res\Ex J,$$
which is the smallest $G$-invariant ideal of $A\times G$ containing $J$,
by \cite[Proposition~11(ii)]{GreLS}.
\end{proof}

%====================================================================

\end{document}